\documentclass[12pt]{article}
\usepackage{epsfig,amssymb,amsmath,amsthm}
\usepackage{graphicx}
\usepackage{verbatim}
\usepackage{hyperref,url}
\newcommand{\nocopyright}{
No Copyright\thanks{
The authors hereby waive all copyright
and related or neighboring rights to this work,
and dedicate it to the public domain.
This applies worldwide.
}}
\title{Conway's doughnuts}
\author{Peter Doyle \and Shikhin Sethi}
\date{Version 1.2.1 dated 11 April 2018\\
\nocopyright
}
\newtheorem{theorem}{Theorem}
\newtheorem{prop}[theorem]{Proposition}
\newtheorem{lemma}[theorem]{Lemma}

\def\RR{\mathbf{R}}

\newcommand{\fig}[3]{
\begin{figure}[htb]
\centerline{
\includegraphics[width=200pt]{figures/#1.pdf}
}
\caption{#3}
\label{#2}
\end{figure}
}
\newcommand{\figsize}[4]{
\begin{figure}[htb]
\centerline{
\includegraphics[width=#1]{figures/#2.pdf}
}
\caption{#4}
\label{#3}
\end{figure}
}

\newcommand{\figbig}[3]{
\begin{figure}[h]
\centerline{
\includegraphics[width=150pt]{figures/#1.pdf}
}
\caption{#3}
\label{#2}
\end{figure}
}

\newcommand{\proofstart}{{\bf Proof.\ }}
\newcommand{\mathproofend}{\quad \qed}
\newcommand{\proofend}{$\quad \qed$}

\begin{document}

\maketitle

\begin{abstract}
Morley's Theorem about angle trisectors
can be viewed as the statement
that a certain diagram `exists',
meaning that triangles of prescribed shapes
meet in a prescribed pattern. 
This diagram is the case $n = 3$
of a class of diagrams we call `Conway's doughnuts'.
These diagrams can be proven to exist using John Smillie's holonomy method,
recently championed by Eric Braude:
`Guess the shapes; check the holonomy.'
For $n = 2, 3, 4$ the existence of the doughnut happens
to be easy to prove because the hole is absent or triangular.
\end{abstract}

\vspace{1em}
\centerline{
%\fbox{
`Watch the doughnut, not the hole.' --- Burl Ives
%}
}

\section{Introduction}
Morley's Theorem about angle trisectors
can be viewed as the statement that a certain diagram `exists',
meaning that triangles of prescribed shapes meet in a prescribed pattern.
(Figure
\ref{fig:nut3}.)
This diagram is the case $n=3$ of a class of diagrams
we call `Conway's doughnuts'.
(This is the only approved spelling;
if you can't abide the `ugh', you may call them `Conway's bagels'.)
See Section 
\ref{sec:flip} for a flip book animation.

Conway's doughnuts can be proven to exist using 
John Smillie's holonomy method:
`Guess the shapes; check the holonomy.'
(See Section
\ref{sec:fluff}.)
Recently this method has been championed by Eric Braude
\cite{braude:morley},
who has used it exactly as we use it here.
For $n=2,3,4$ the existence of the doughnut happens to be
easy to prove because the hole is absent or triangular.

\section{Bisectors} \label{sec:two}

The Bisector Theorem
(Euclid's Elements, Book IV, Proposition 4)
says that the three
angle bisectors of a Euclidean triangle
pass through a common point,
namely the incenter of the triangle.
This is easy to understand and prove:
$x=y$ and $y=z$ implies $x=z$.
With a view to generalizations, we rephrase it as follows.

Let $\Delta(A,B,C)$ denote the similarity class of Euclidean triangles
having angles $A,B,C$ in counter-clockwise order.
Say that a diagram like those
in Figures
\ref{fig:nut2} and \ref{fig:nut3}
\emph{exists}
if triangles of the designated shapes
fit together as indicated.

\begin{theorem}[Bisector Theorem]
Let $x' = x + \tau/4$ for any $x \in \RR$.
Assume $a+b+c=0'=\tau/4$.
The diagram of Figure
\ref{fig:nut2}
exists,
in the sense that triangles of shapes
$\Delta(a,b,c'),\Delta(b,c,a'),\Delta(c,a,b')$
fit together around a vertex,
forming a triangle of shape $\Delta(2a,2b,2c)$.
\end{theorem}

\fig{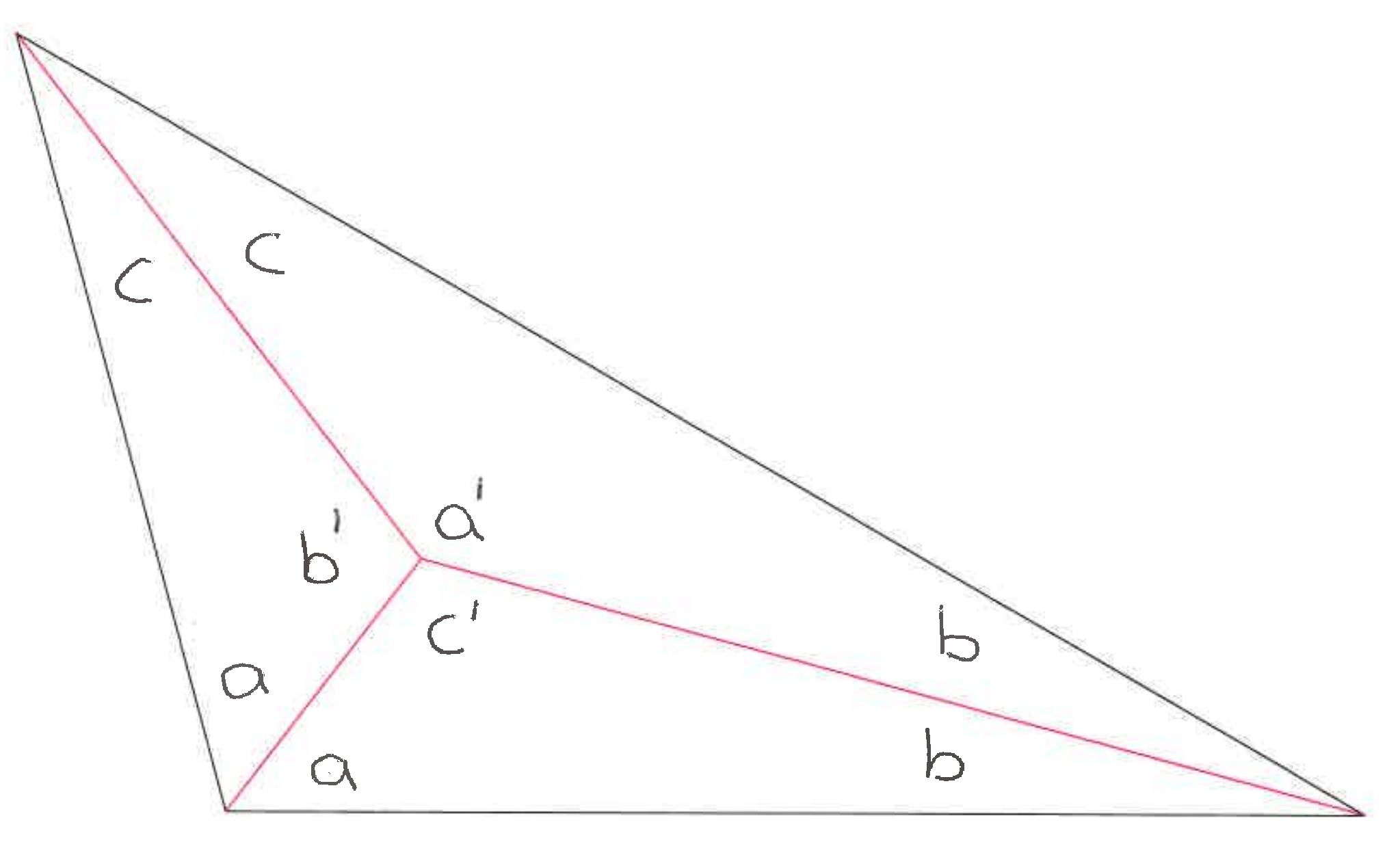}{fig:nut2}{
The Bisector Theorem.
$x'=x+\tau/4$, $a+b+c=0'$.
}

\proofstart
We use Smillie's holonomy method:
`Guess the shapes; check the holonomy.'
The shapes of triangles we've specified are legal because
$a+b+c'=0''=\tau/2$, etc.
Do they fit together as advertised?
We check the holonomy around the interior vertex.
If we nail down some instance $\delta_0$ of one of the subtriangles,
say $\Delta(a,b,c')$,
that nails down an instance $\delta_1$ of its neighbor $\Delta(b,c,a')$,
which nails down an instance $\delta_2$ of its neighbor $\Delta(c,a,b')$,
which nails down an instance $\delta_3$ of its neighbor $\Delta(a,b,c')$.
We need $\delta_3 = \delta_0$.
Their $c'$ vertices match, and their orientations match because
the angle sum around the vertex is $a'+b'+c'=0''''=\tau$.
Do their sizes match?
Yes,
because by the law of sines
the ratio of the $(a,c')$ side of $\delta_3$
to that of $\delta_0$ is
\[
\frac{\sin(a)}{\sin(b)}
\frac{\sin(b)}{\sin(c)}
\frac{\sin(c)}{\sin(a)}
=1
.
\mathproofend
\]

Let us package the computation done here.
(Cf.\ Braude \cite[section 2.3]{braude:morley}).)
\begin{lemma}[EZ Holonomy]
The holonomy of a collection of shapes
\[
\Delta(a_1,b_1,c_1),\ldots,\Delta(a_n,b_n,c_n)
\]
arranged in a ring around their $c_i$ vertices is trivial if
$c_1+\ldots+c_n = \tau$ and the `trailing' angles $a_1,\ldots,a_n$
are a permutation of the `leading' angles $b_1,\ldots,b_n$.
\proofend
\end{lemma}

\section{Trisectors} \label{sec:three}

Moving from bisectors to trisectors we have Morley's Theorem
\cite{morley:trisectors}.

\begin{theorem}[Morley's Theorem]
The diagram of 
Figure
\ref{fig:nut3}
exists.
\end{theorem}

\fig{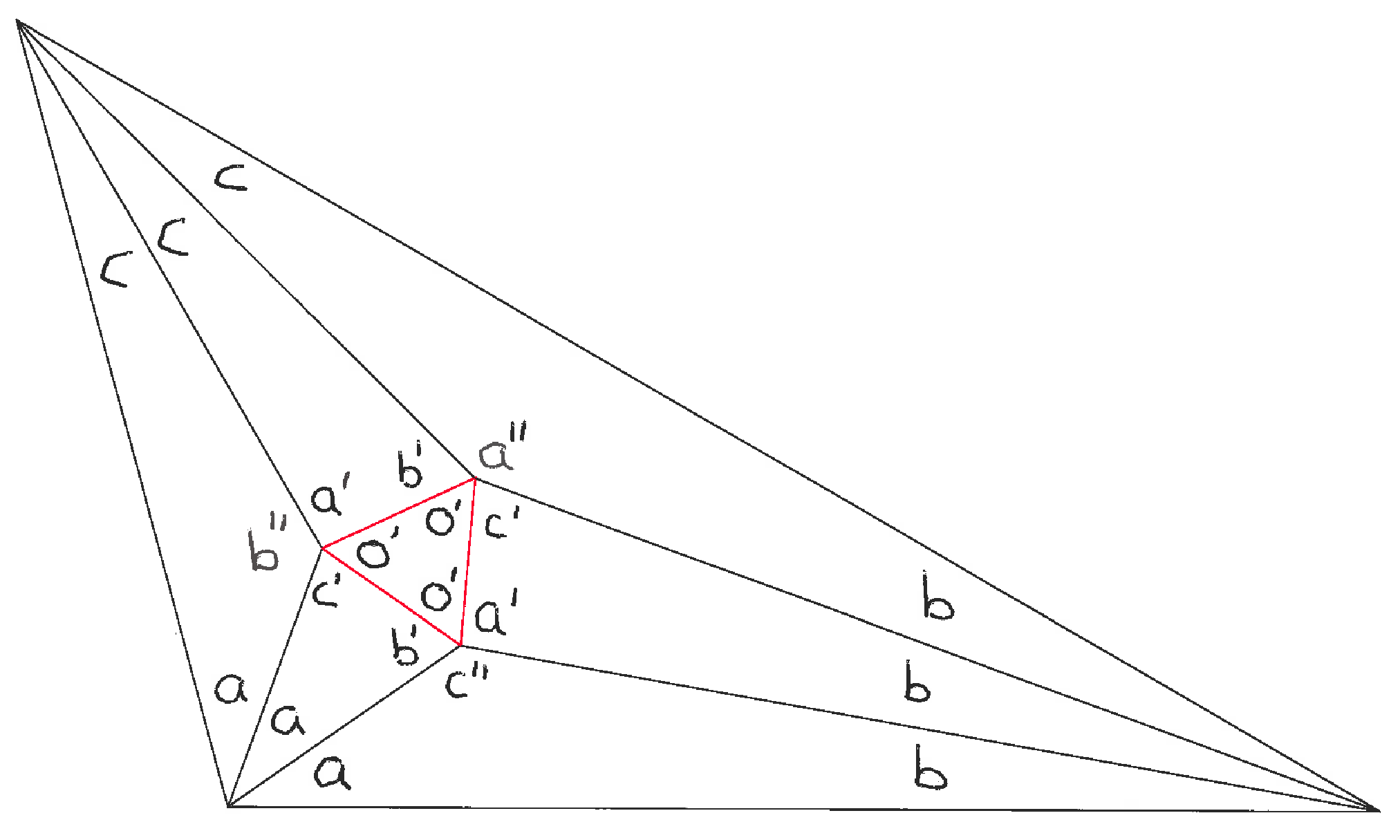}{fig:nut3}{
Morley's Theorem.
$x'=\tau/6$, $a+b+c=0'$.
}

\proofstart
The angles add to $\tau/2$ in all the triangles,
and the EZ Holonomy condition holds at the interior vertices
(by symmetry we need only check one vertex).
\proofend

\section{Quadrisectors} \label{sec:four}

Bill Thurston observed that topologically,
the diagrams for the Bisector Theorem
and Morley's Theorem are the planar maps of a tetrahedron and an
icosehedron,
and wondered (though not enough to think about the question himself)
what the analog for an icosahedron might be.
John Conway promptly supplied
Figure
\ref{fig:nut4icos},
though with some misgivings since it is really just Figure
\ref{fig:nut4}
with some extraneous lines.
Figure
\ref{fig:nut4icosvar} shows a third related configuration.
All these diagrams exist by EZ Holonomy.
We take 
Figure 
\ref{fig:nut4}
as the proper extrapolation to quadrisectors,
and call it Conway's Theorem.

\figbig{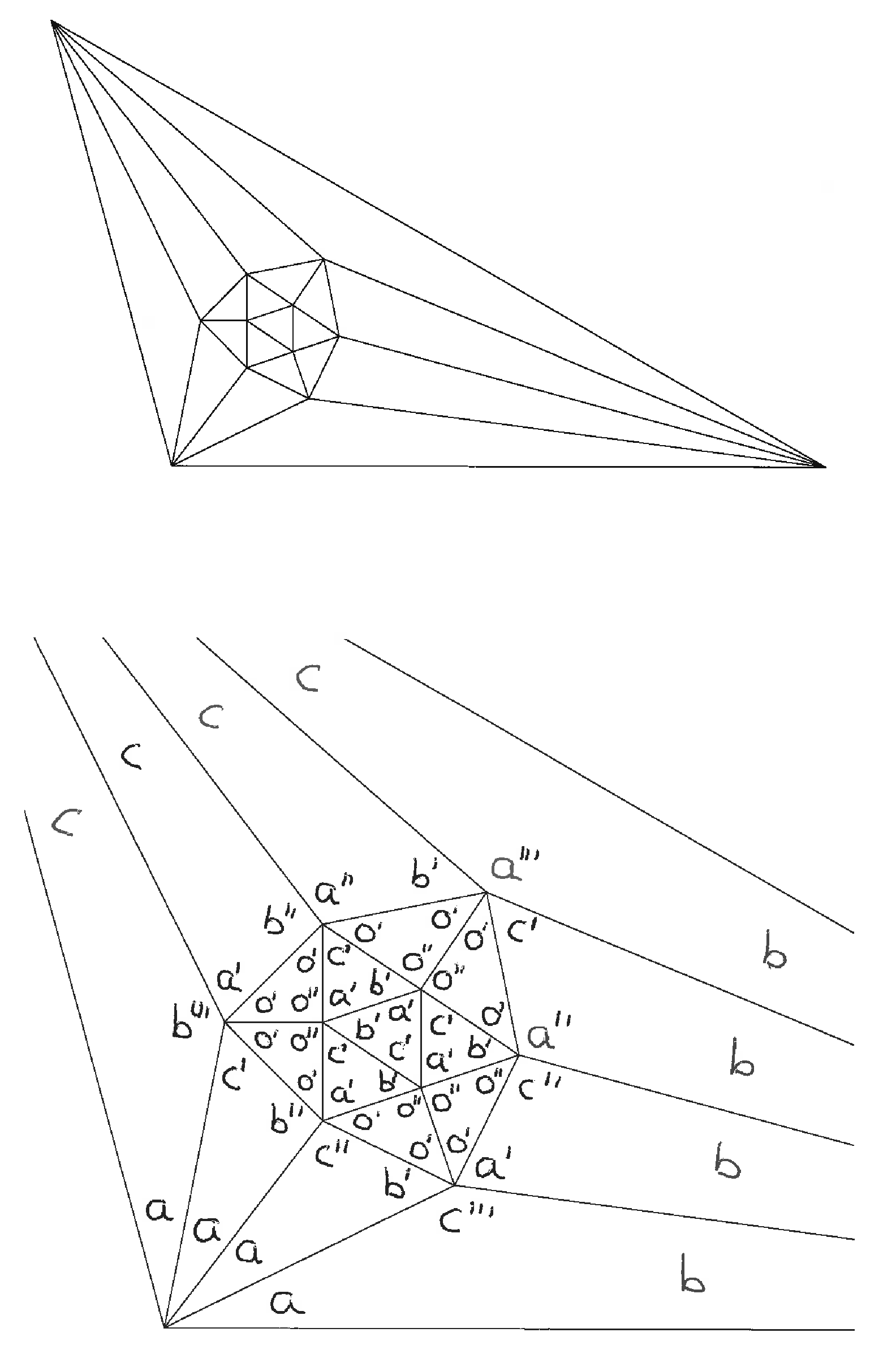}{fig:nut4icos}{
Conway's icosahedron.
$x'=\tau/8$, $a+b+c=0'$.
}
\figbig{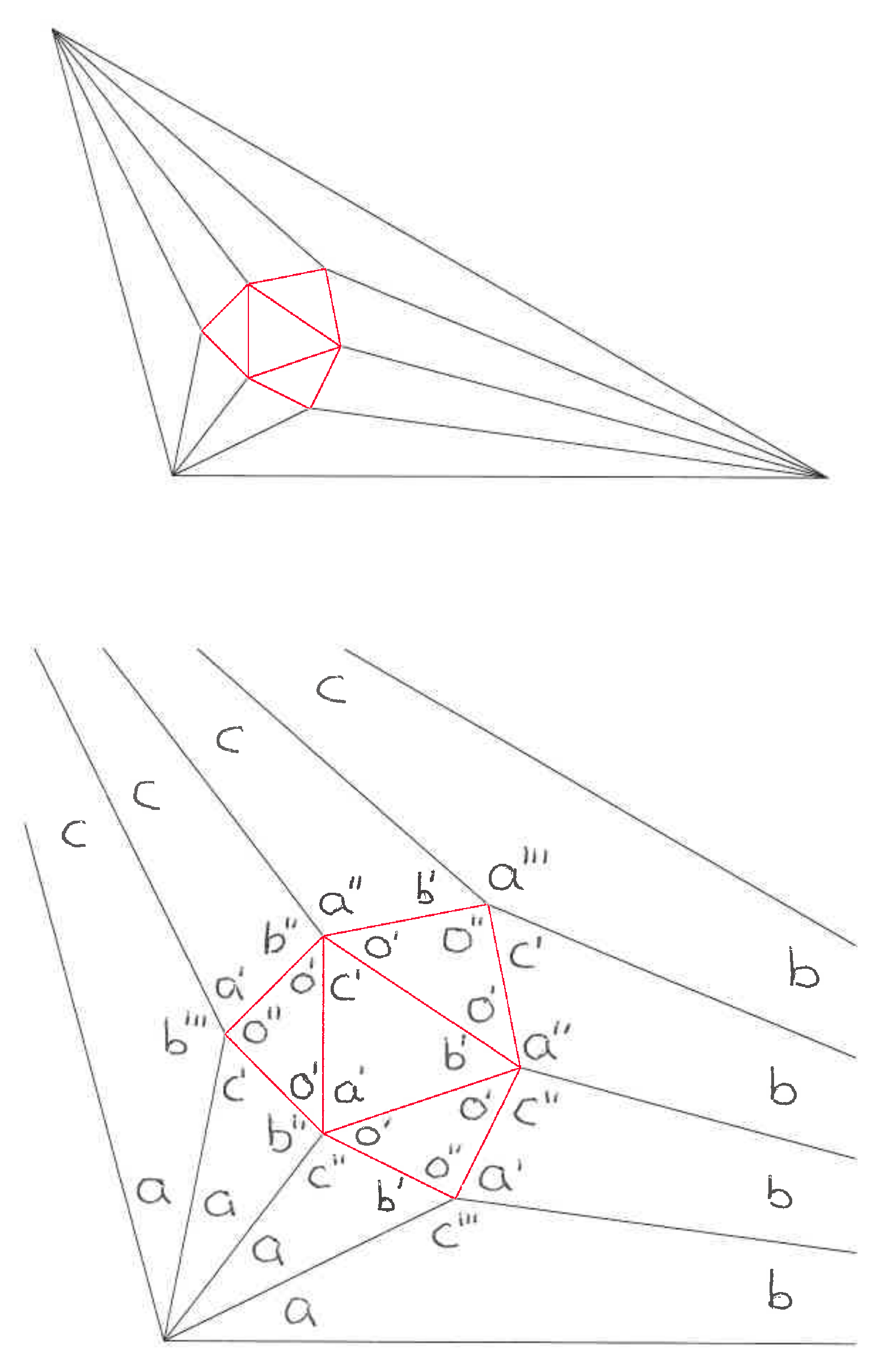}{fig:nut4}{
Conway's Theorem.
$x'=\tau/8$, $a+b+c=0'$.
}
\figsize{200pt}{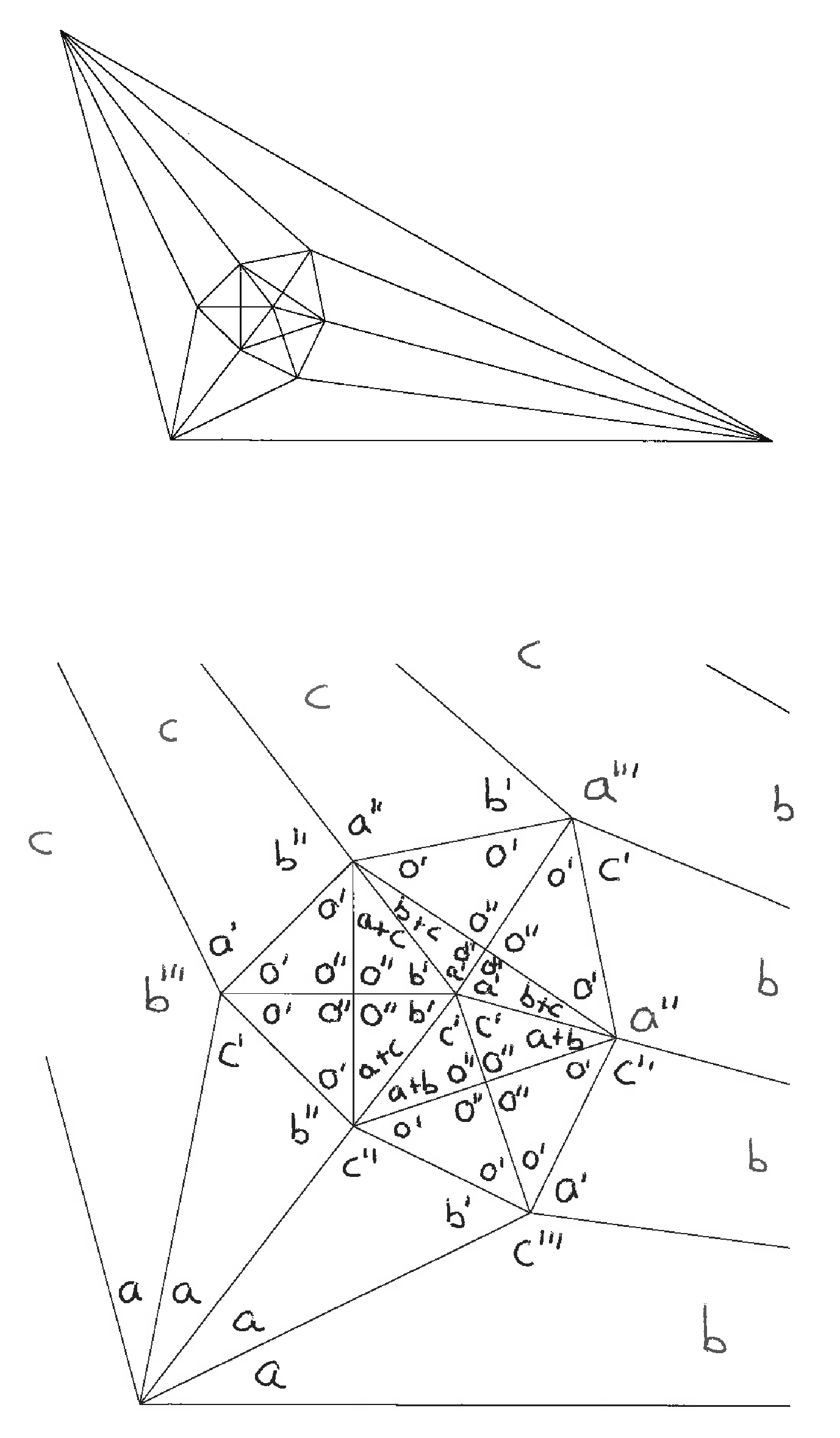}{fig:nut4icosvar}{
A variation on Conway's icosahedron.
$x'=\tau/8$, $a+b+c=0'$.
This diagram is a subdivision of the angle bisector diagram.
The circumcenter of the three new points on the angle bisectors
coincides with the incenter of the triangle.
}
\begin{theorem}[Conway's Theorem]
The diagram of Figure
\ref{fig:nut4}
exists.
\end{theorem}

\proofstart
EZ Holonomy.
\proofend

\clearpage

\section{Conway's doughnuts} \label{sec:five}
What about
$n=5$?
(Let's dispense with the word `quinquisector'.)
Thanks to Conway, it's clear what to put around the outside.
(Figure
\ref{fig:nut5outside}.)
\figbig{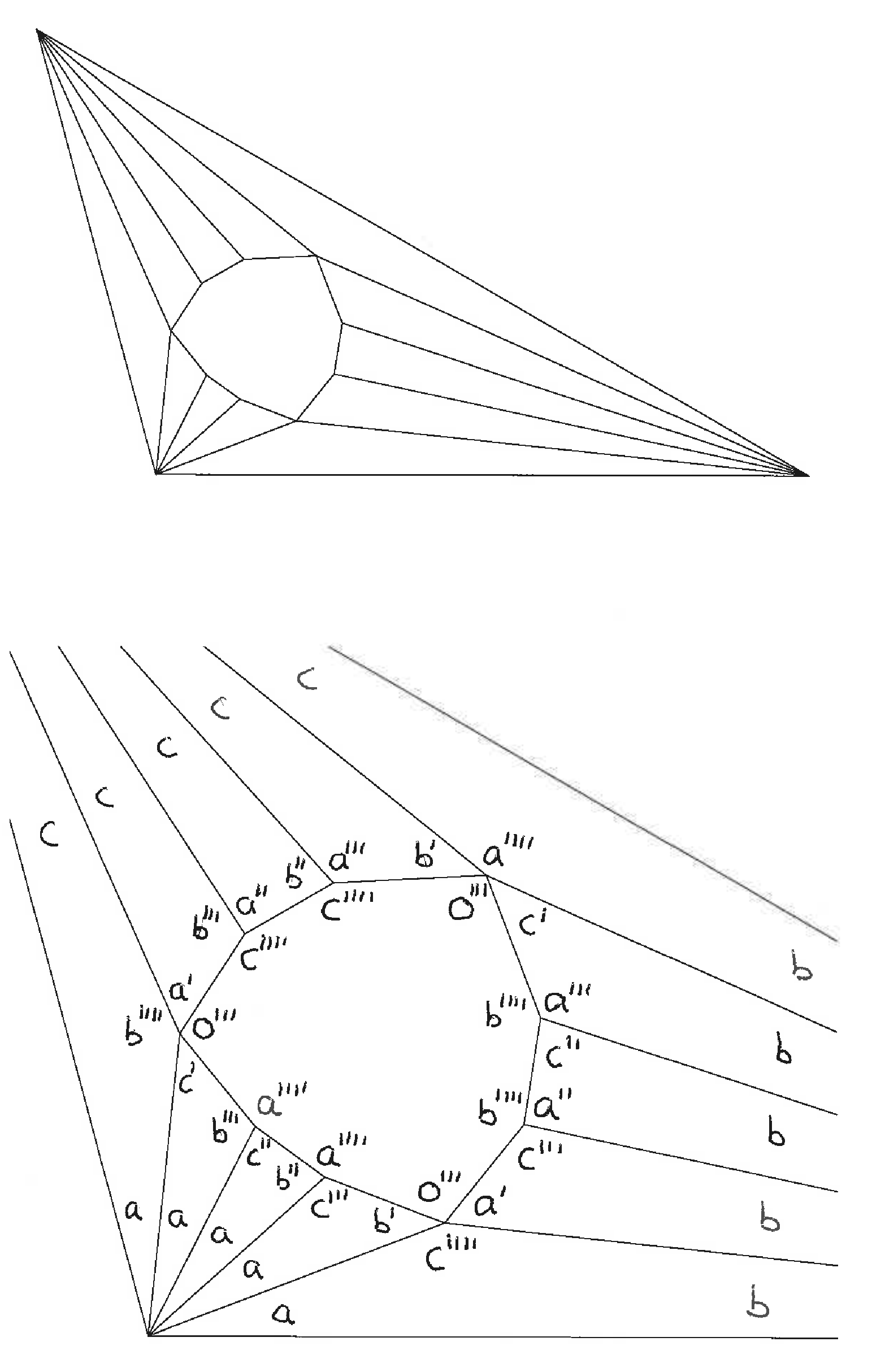}{fig:nut5outside}{For $n=5$, it's clear what should go outside.}
How to fill in the hole?
Beyond adding isosceles triangles as in
Figure
\ref{fig:nut5},
\figbig{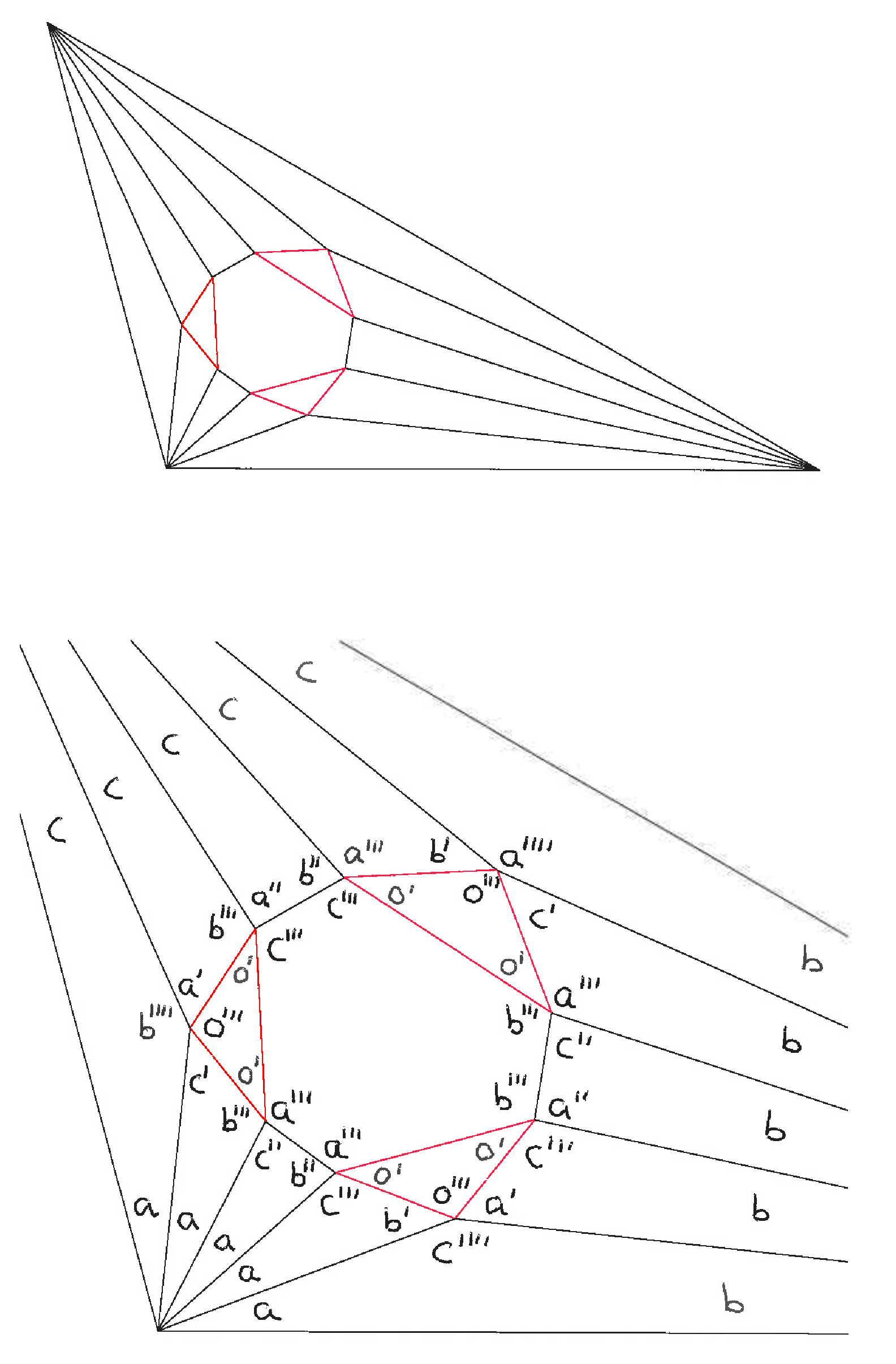}{fig:nut5}{Conway's 5-doughnut.}
there seems to be nothing good.
Let's face it, we've run out of Platonic solids.
(We did try to do something with the tiling of the Euclidean
plan by triangles meeting six to a vertex, but couldn't
make this work.)

So let's declare a victory.
Forget the hole.
Watch the doughnut.

\begin{theorem}[Conway's doughnuts]
The analog of Figure
\ref{fig:nut5}
exists for all $n$.
\end{theorem}

\proofstart
This comes down to checking that the subtriangles meeting at the corners
fit properly.
We isolate this in the next Proposition.
\proofend

\begin{prop}
Suppose $a+b+c=\tau/(2n)$.
Triangle shapes
\[
\Delta(a,b,c+\frac{n-1}{2n}\tau),
\Delta(a,b+\frac{1}{2n}\tau,c+\frac{n-2}{2n}\tau),
\ldots,
\Delta(a,b+\frac{n-1}{2n}\tau,c)
\]
\begin{comment}
\[
\Delta(a,b,c+(n-1)\frac{\tau}{2n}),
\Delta(a,b+\frac{\tau}{2n},c+(n-2)\frac{\tau}{2n}),
\ldots,
\Delta(a,b+(n-1)\frac{\tau}{2n},c)
\]
\end{comment}
fit properly inside
$\Delta(na,nb,nc)$
as illustrated in Figure
\ref{fig:nut5corner}
for the case $n=5$.
\end{prop}

\figsize{370pt}{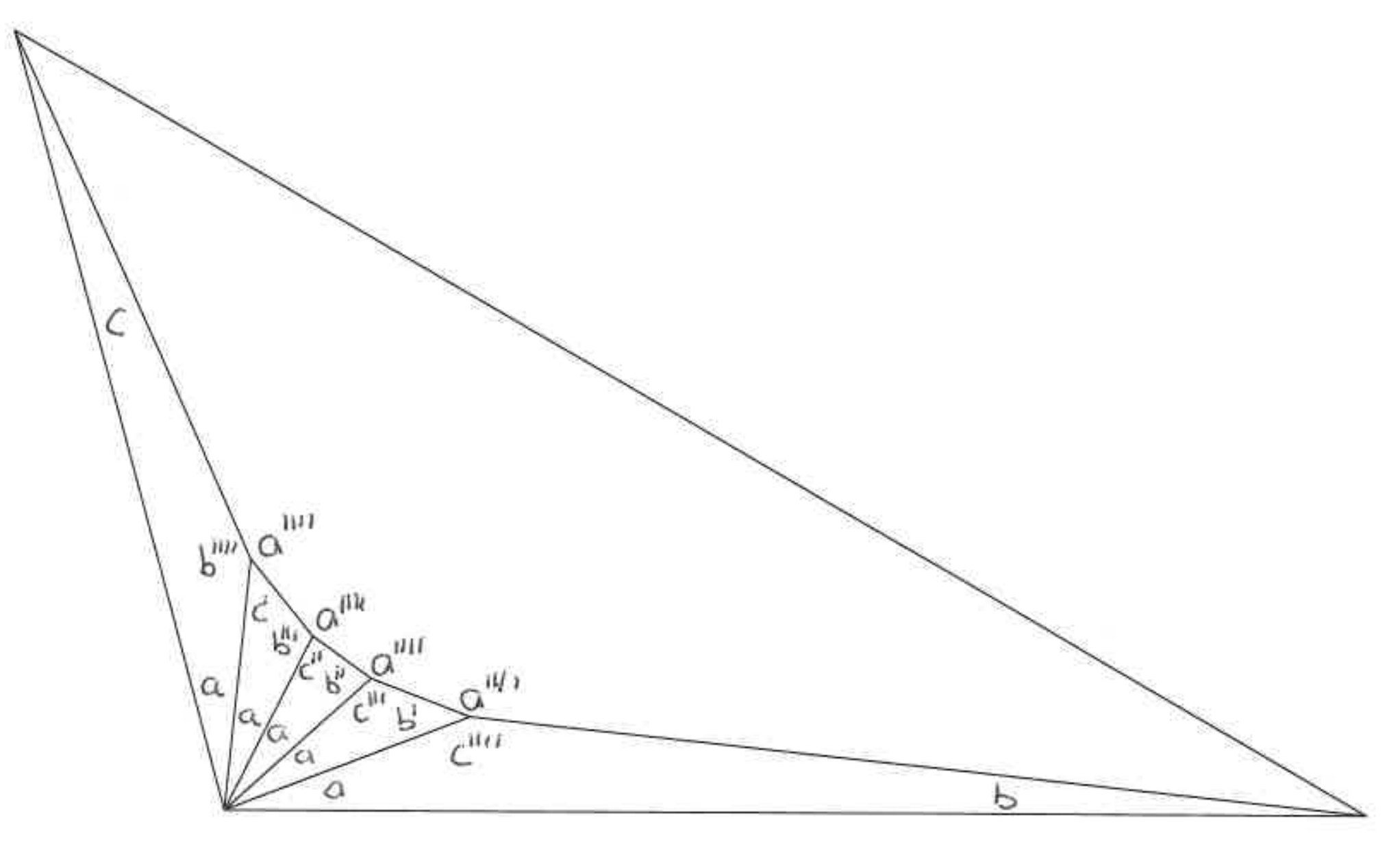}{fig:nut5corner}{We just need to show that this corner
figure exists.}

\proofstart
We take the case $n=3$ as a model.
We're trying to prove the existence of the diagram
consisting of our three triangles
$\Delta(a,b,c''),\Delta(a,b',c'),\Delta(a,b'',c')$,
together with the large $\Delta(3a,3b,3c)$,
glued together to form a cone flattened down into the plane.
The holonomy condition at the vertex is that
\[
\frac{\sin(3b)}{\sin(3c)}
=
\frac{\sin(b)}{\sin(c'')}
\frac{\sin(b')}{\sin(c')}
\frac{\sin(b'')}{\sin(c)}
.
\]
Since we can vary $c$ while fixing $b$, we had better have the identity
\begin{align*}
\sin(3x)
&= A \sin(x) \sin(x') \sin(x'')
\\
&= A \sin(x) \sin(x + \frac{\tau}{6}) \sin(x + 2 \frac{\tau}{6})
,
\end{align*}
for some constant $A$.
This follows from general principles because the two sides are entire functions
of modest growth with the same roots.
\proofend

The trigonometric identity that we've proven analytically
up to a constant factor is the following identity, which is key to
computations of hyperbolic volume.
(See Milnor
\cite{milnor:hyperbolic}.)

\begin{prop}[Cyclotomic identity]
\[
\sin(n \theta) =
2^{n-1}
\sin(\theta) 
\sin(\theta+\frac{1}{2n}\tau)
\ldots
\sin(\theta+\frac{n-1}{2n}\tau)
.
\]
\end{prop}

\proofstart
Again we take the case $n=3$ as a model.
Start with the identity
\[
z^3-1 = (z-1)(z-\omega)(z-\omega^2)
,
\]
with $\omega = e^{i \tau/3}$.
Plug in $z=e^{i 2 \theta}$, take absolute values,
and use the fact that for $\alpha, \beta \in \RR$
\[
|e^{i \alpha} - e^{i \beta}| =
2 |\sin((\alpha-\beta)/2)|
\]
to get
\begin{align*}
2 |\sin(3 \theta)|
&=
2 |\sin(\theta)|
\cdot
2 |\sin(\theta-\tau/6)|
\cdot
2 |\sin(\theta-2\tau/6|
\\&=
2 |\sin(\theta)|
\cdot
2 |\sin(\theta+2\tau/6)|
\cdot
2 |\sin(\theta+\tau/6)|
.
\end{align*}
Remove absolute values and check the sign to get the desired identity.
\proofend

This identity cries out for a geometric proof---see 
Section
\ref{sec:isthatall}
below.

\clearpage
\section{Asymptotics}

Properly renormalized, 
the polygonal boundary of the doughnut hole
near a vertex of angle $A$ approaches
the curve $\mathrm{Im}(z^{\tau/(2A)})=1.$
For a right angle $A = \tau/4$ this limiting curve is the
hyperbola $2xy=1$;
for a 60-degree angle $A=\tau/6$ it's the `monkey hyperbola'
$3x^2 y - y^3=1$.
In general it's a curve that looks like a hyperbola by isn't.
To demonstrate, in
Figure \ref{fig:nut10bent}
we've raised one of our doughnuts to a fractional power to straighten
out one of the corners of the triangle.
(See also the flip book in Section
\ref{sec:flip}.)

\figsize{370pt}{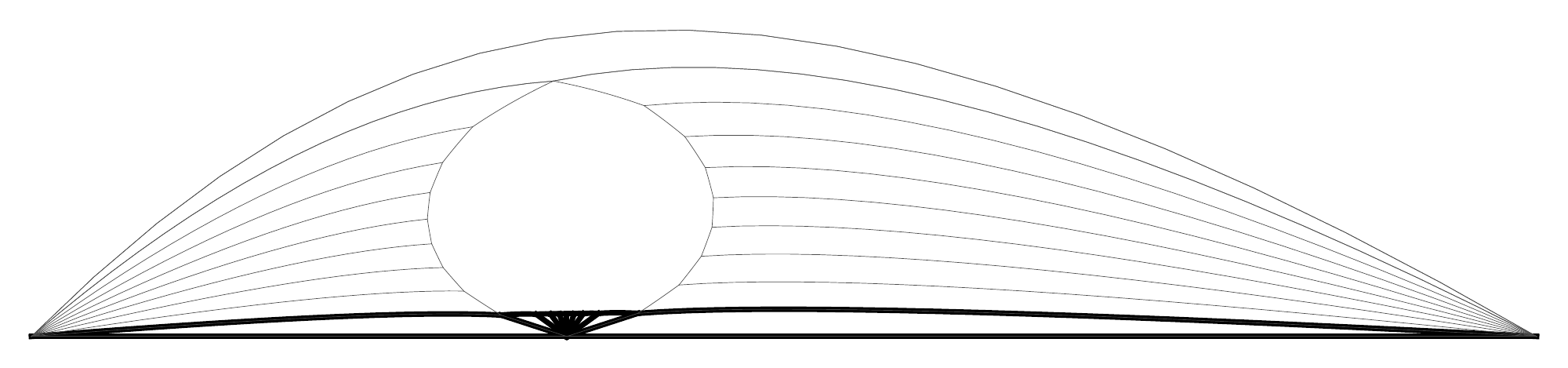}{fig:nut10bent}
{A 10-doughnut raised to a fractional power.}

\section{Reflections by Peter Doyle} \label{sec:fluff}

The holonomy approach to proving geometrical theorems
is to `guess the shapes; check the holonomy.'
In 1976 John Smillie introduced this method to construct the developing map
of an affine torus obtained by glueing opposite sides of an arbitrary
Euclidean quadrilateral.
Bill Thurston liked the picture and included it in his 1980 notes
on 3-manifolds
\cite{thurston:three}.
(Figure
\ref{fig:wptsmillie}.)
I found this idea striking, 
and took note of it, though I missed the fact
(stated back on page 3.3) that the idea was due to John Smillie.
It was only very recently that I looked back at Bill's notes,
and saw
that it was John Smillie rather than Bill who pioneered the holonomy
method.

\figsize{370pt}{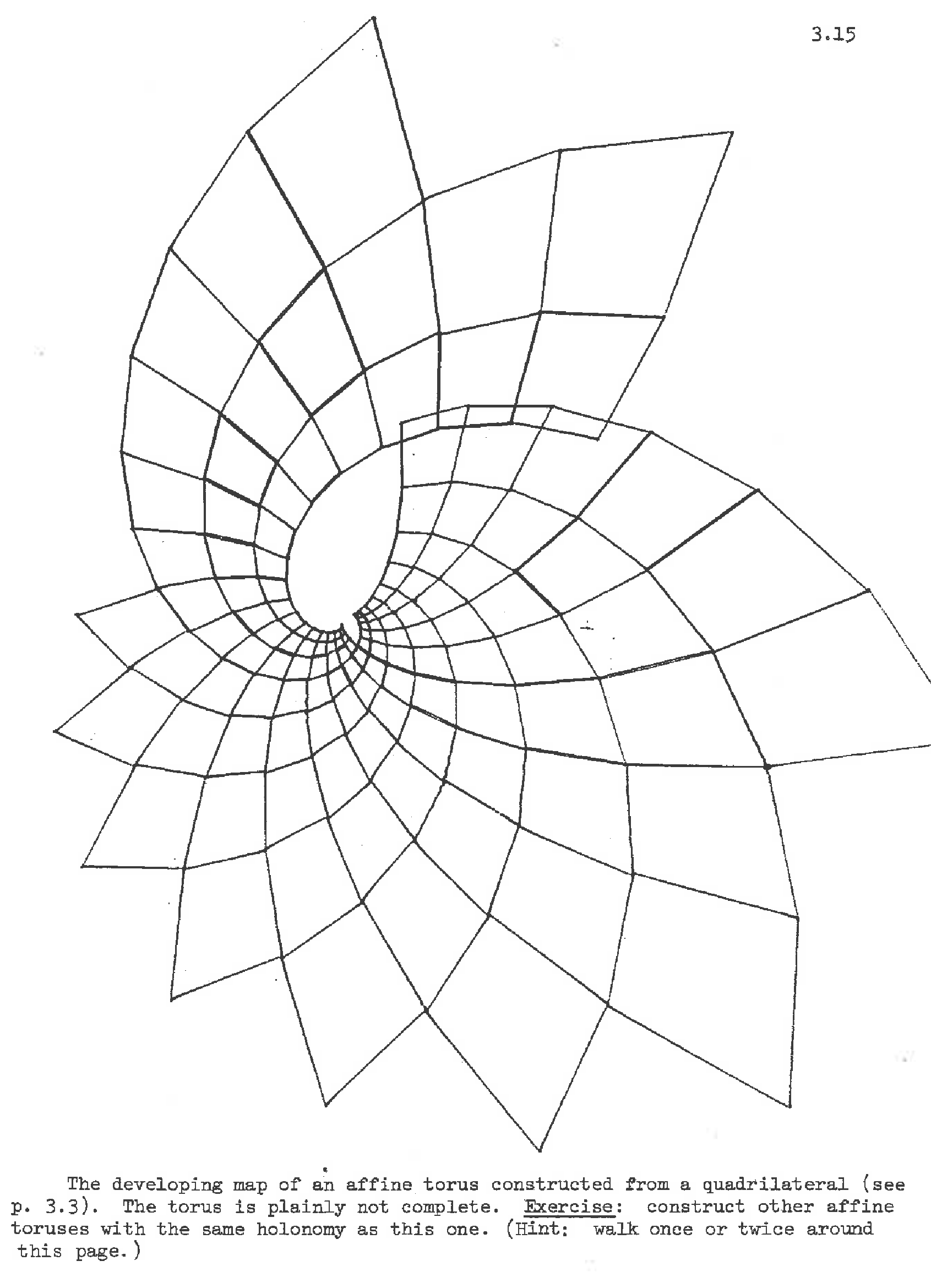}{fig:wptsmillie}{John Smillie's construction of affine tori,
from Thurston
\cite{thurston:three}.}

Some years later when circle packings were all the rage
I used this holonomy method to prove the existence of 
`exponential circle packings',
which are the analog for circle packings of the exponential function.
(See Figures
\ref{fig:exppack},\ref{fig:expangles},\ref{fig:expfit}.)
%Someone decided to call special cases of these packing `Doyle spirals'.
I presumed (and still presume) that Bill already knew about these packings,
though I never asked him.

\figsize{400pt}{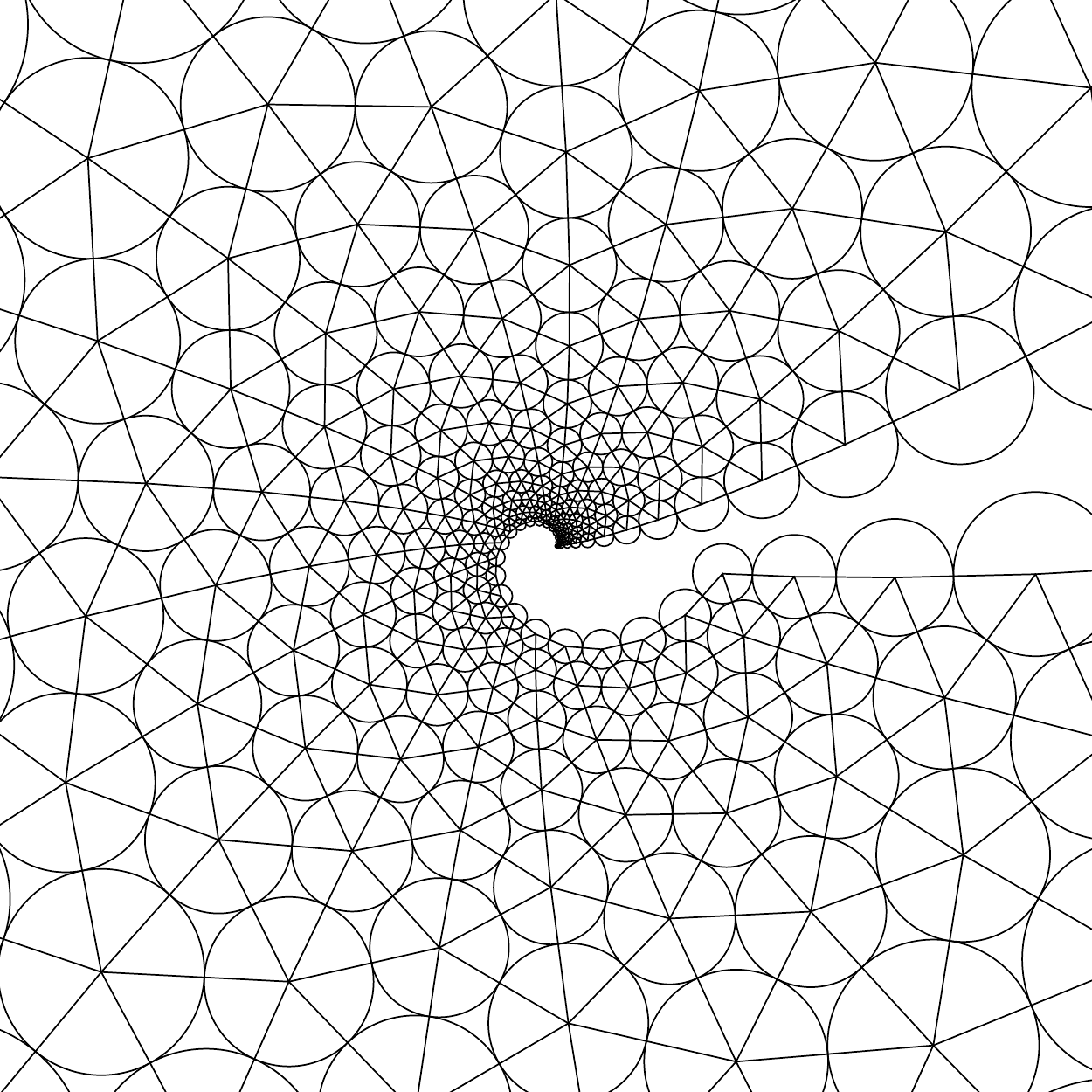}{fig:exppack}{An exponential circle packing.}

\figsize{120pt}{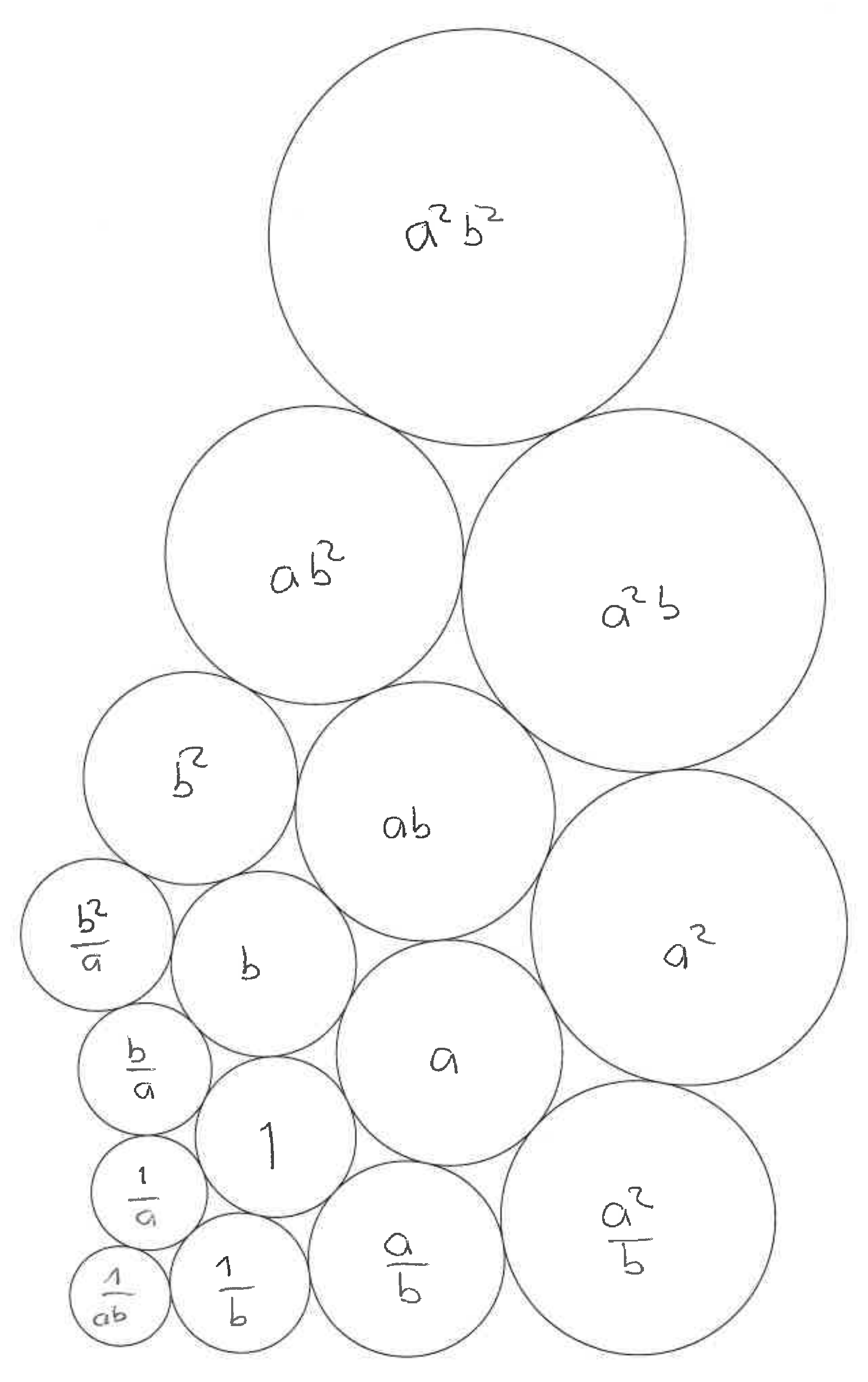}{fig:exprad}{Radii of circles in an exponential packing.}

\figsize{110pt}{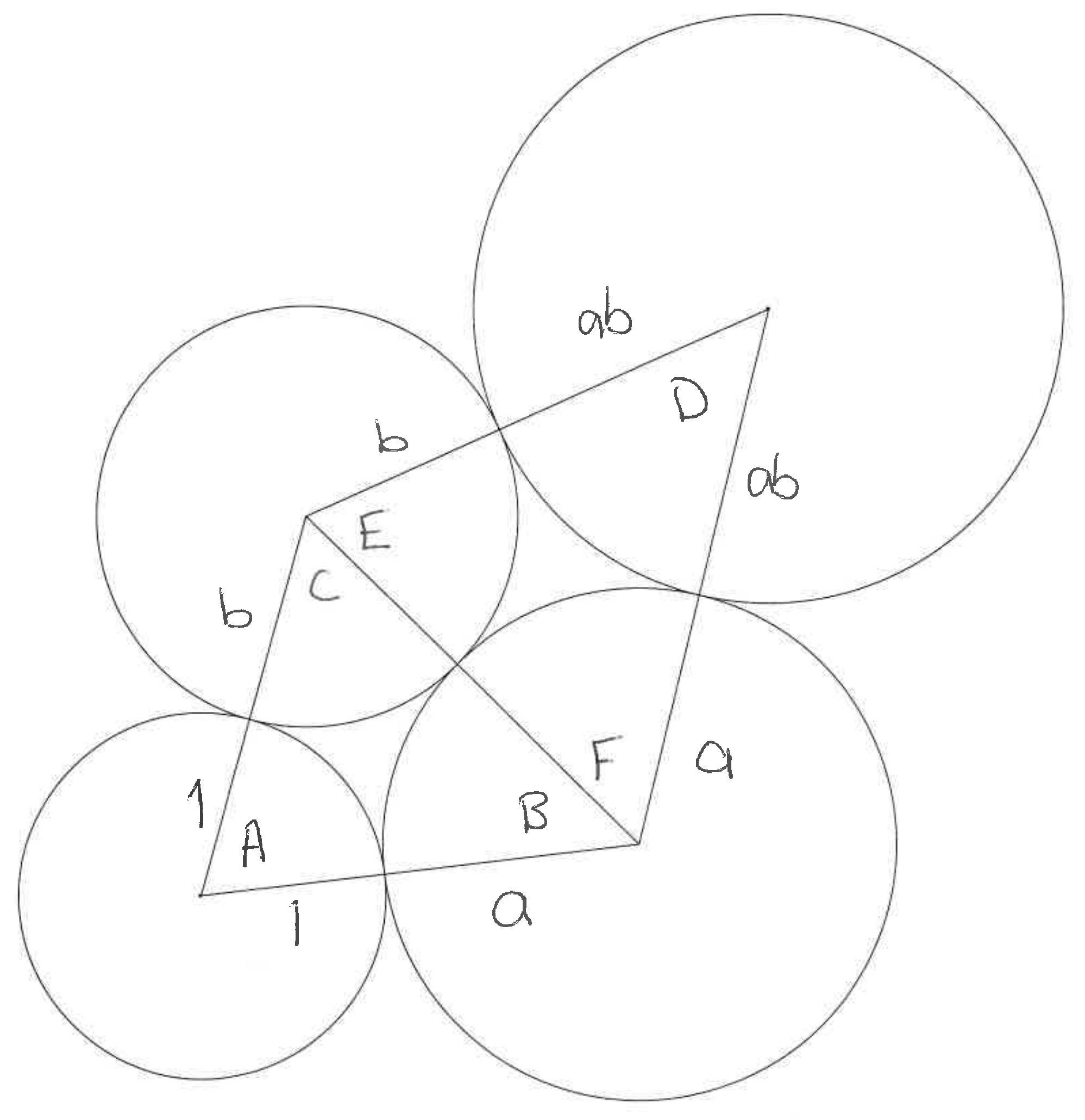}{fig:expangles}{Shapes of triangles in an exponential packing.}

\figsize{110pt}{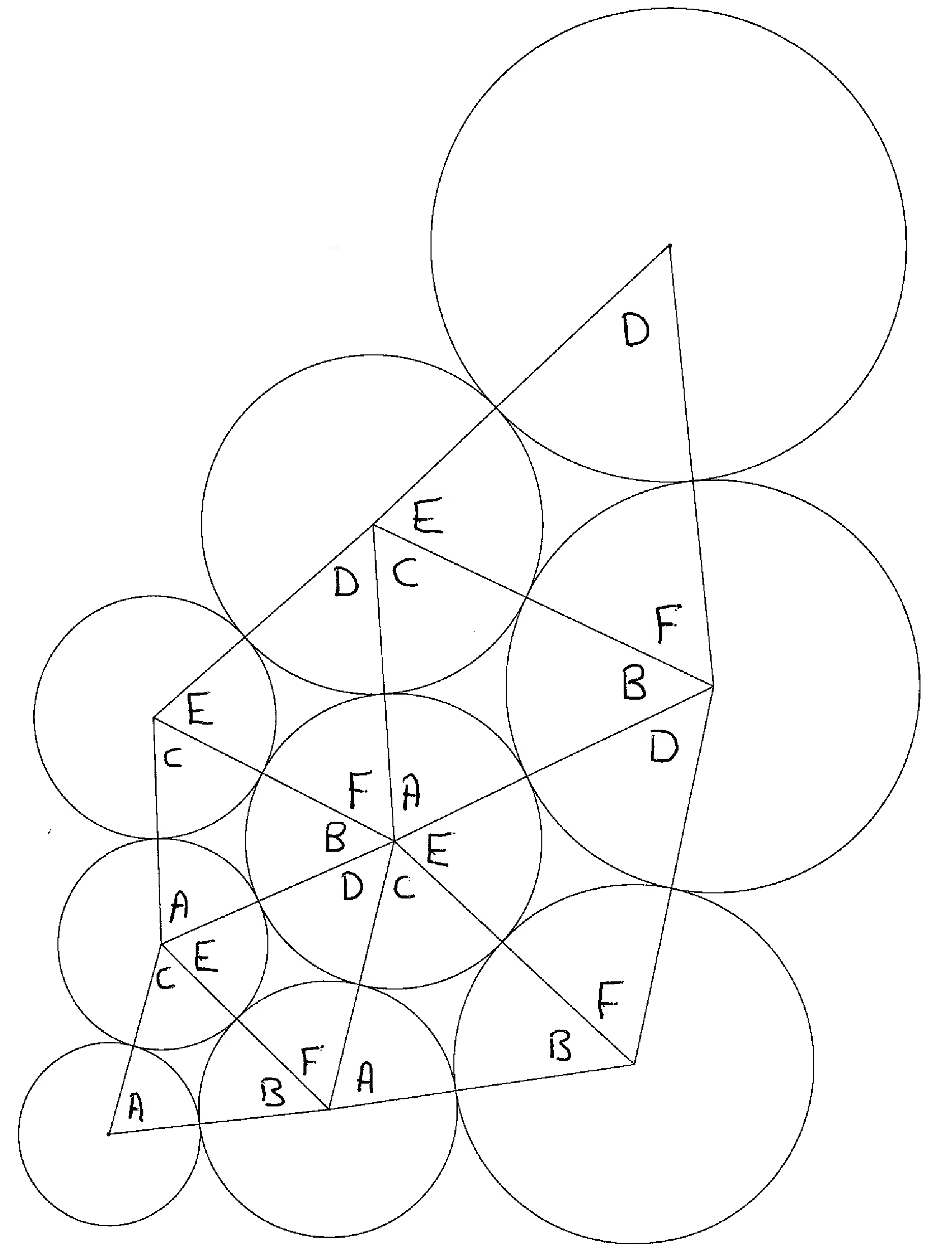}{fig:expfit}{
Triangles fit around a vertex by EZ holonomy.
This is true for any two triangle shapes, not just these special shapes.
This brings us right back to Smillie's
affine tori.
In the space of affine tori, we're dealing with a subvariety
consisting of perturbations of a hexagonal Euclidean torus.
}

\clearpage

In 1990/1991 when Bill Thurston, John Conway, Jane Gilman, and I developed
the course `Geometry and the Imagination'
\cite{gi:mpls}
 at Princeton,
we took to starting each class with a short presentation by John of
one of `Conway's Fascinating Facts'.
Among these was Morley's Theorem.
I recognized that Morley's Theorem 
could be proved using the holonomy method
(which I still thought was due to Bill),
as could other standard facts about triangles.
(See Figures
\ref{fig:centers}, \ref{fig:isogonal}, \ref{fig:chopsticks}.)
John Conway made the holonomy proof more elementary
by adding extra lines to the Morley diagram.
(See
Conway
\cite{conway:three}.)
I felt these extra lines made the proof more complicated,
and hard to remember:
Without them you can recostruct the proof as long as you remember to
`guess the shapes; check the holonomy'.
It has been gratifying to see subsequent authors strip the extra lines out.
(Cf. Braude
\cite{braude:morley}.)

\figsize{150pt}{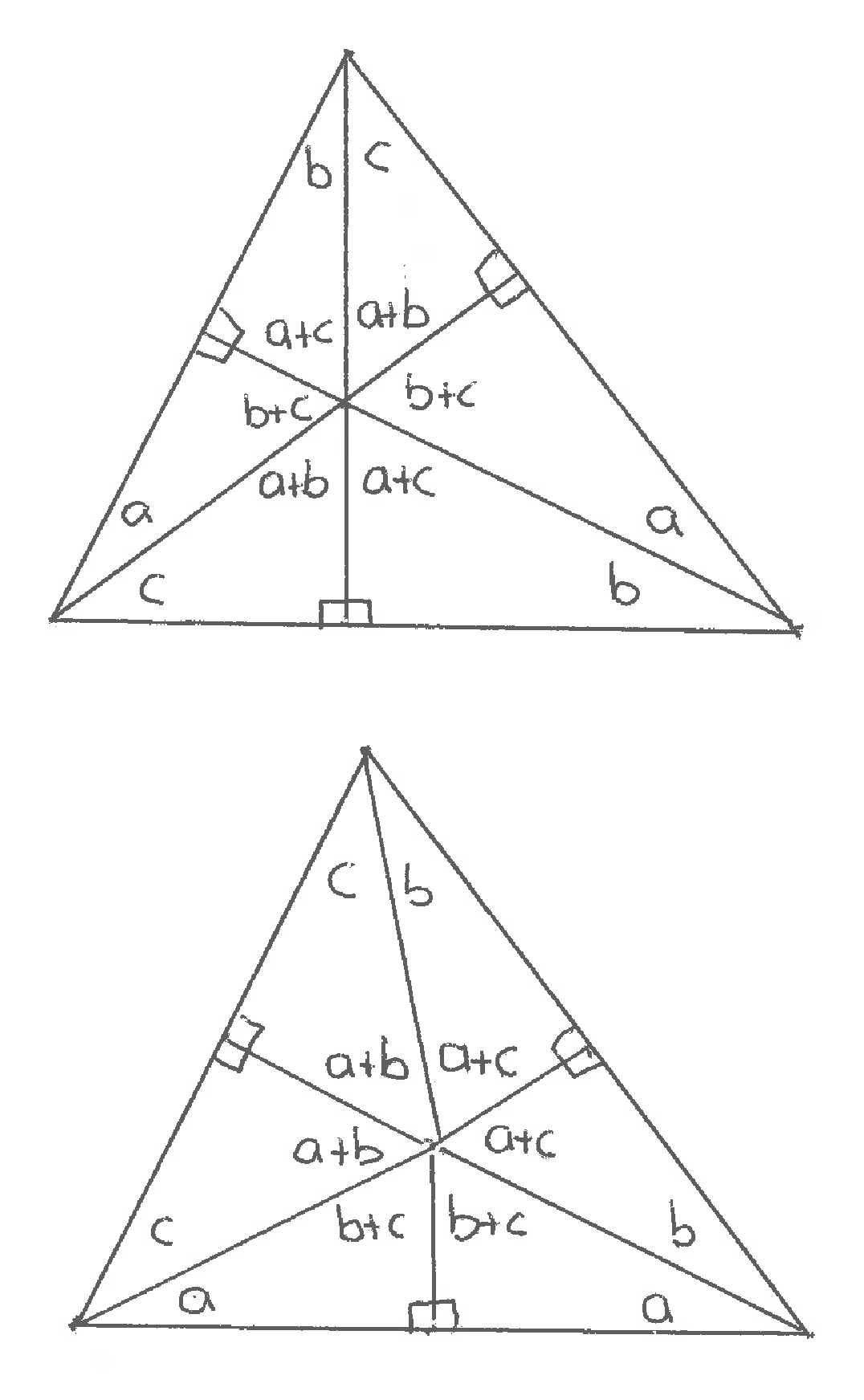}{fig:centers}{
Orthocenter  and Circumcenter.
$a+b+c=\tau/4$.
Note that the second diagram implies the Inscribed Angle Theorem:
The visual angle of an edge as seen from the circumcenter is twice that
as seen from the opposite vertex.
(See Rich Schwartz's beautiful proof
\cite{schwartz:surfaces}
of what he calls
the `X Theorem'.)
}
\figsize{350pt}{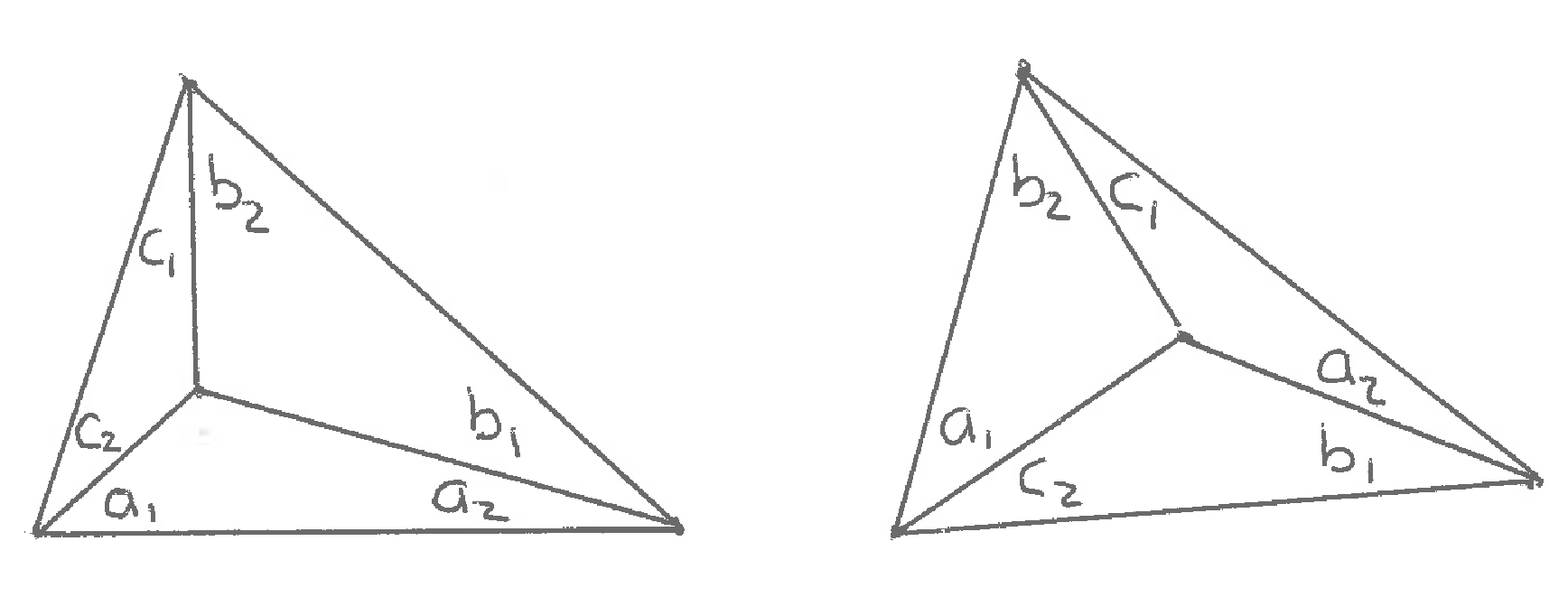}{fig:isogonal}{
Isogonal Conjugates.  If the first diagram exists (with appropriate central
angles), then so does the second.
}
\figsize{350pt}{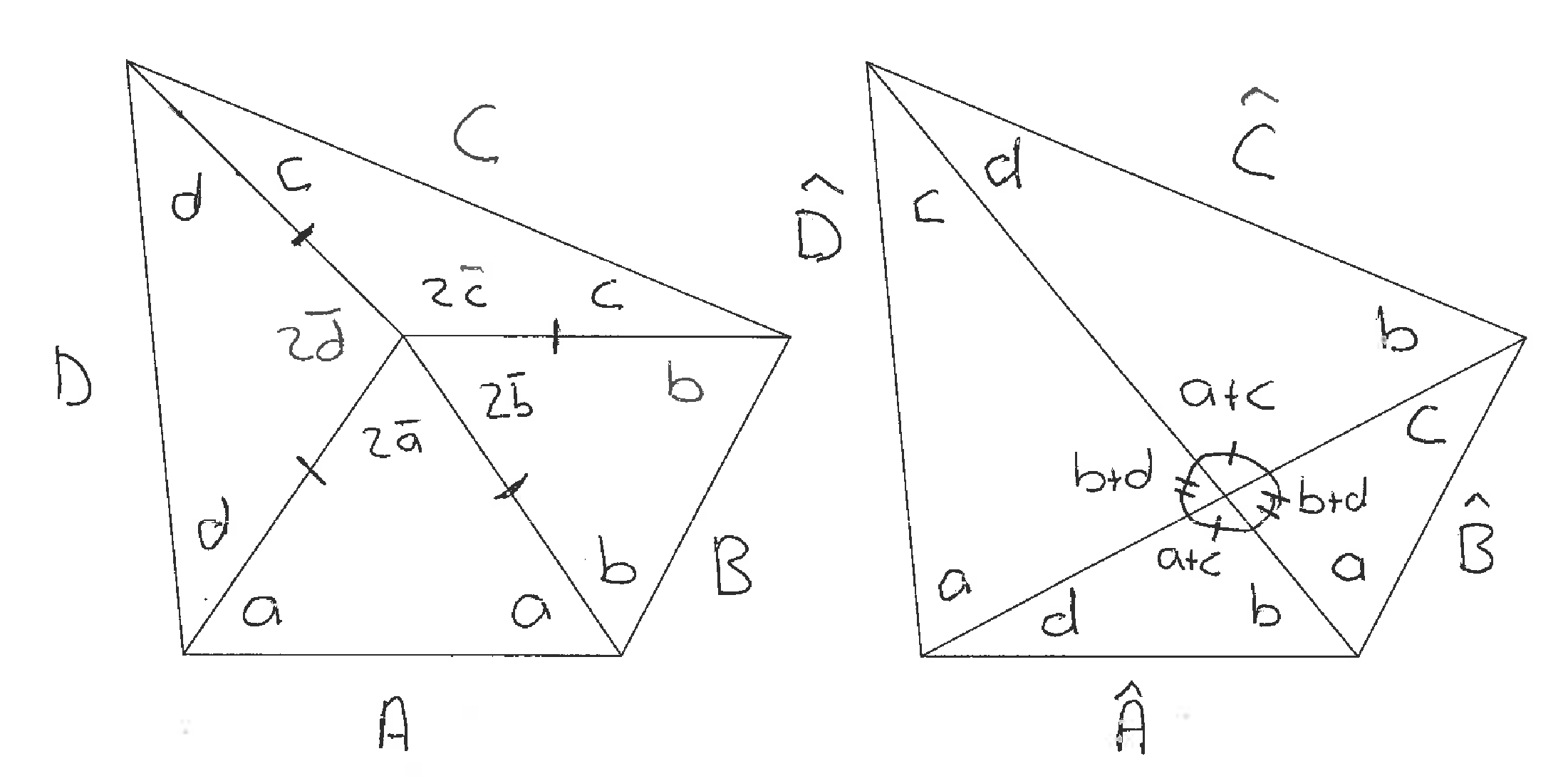}{fig:chopsticks} {Cyclic Quadrilateral.
$\bar{x}=\tau/4-x$, $a+b+c+d=\tau/2$.
The diagram here is a doubled quadrilateral,
with $A$ glued to $\hat{A}$, etc. 
Front and back exist individually by EZ Holonomy.
The front is inscribable in a circle.
Front and back fit together because $B/A = \hat{B}/\hat{A} = \sin(b)/\sin(a)$,
etc.
The back shows that an edge has the same visual angle from either of the two opposite vertices.
This is another aspect of the Inscribed Angle Theorem.
}

\clearpage

As described in Section
\ref{sec:four},
John then produced Conway's Theorem,
aided by the felicitous $a,a',a''$ notation that we have
appropriated here.
I spent decades trying to take the next step to $n=5$,
replacing the icosahedron with tiling of the plane by triangles
meeting six to a vertex.
Finally my coauthor Shikhin Sethi
stepped in and did a computer search, which convinced me that
there is no useful way to add triangles with angles depending linearly
on $a,b,c$.
This led to the thesis proposed here,
that the Bisector
Theorem, Morley's Theorem, and Conway's Theorem are best viewed
as special cases where the hole in Conway's doughnuts disappears
($n=2,3$) or is a triangle ($n=4$).

As you will have noted,
my only contribution to this discussion
was to 
recognize Morley's Theorem as a natural for
the holonomy method,
and as for that, the literature on Morley's theorem is littered
with proofs that are arguably equivalent.
The rest was the work of John Conway.
It was John who (with a hint from Bill)
saw how to push past Morley's Theorem
to Conway's Theorem.
After that the step to Conway's doughnuts is clear;
if John never brought them to market, it was likely because he
felt that without something to fill the holes, they
weren't worth anything.
That's the notion we're trying to combat here.

\section{Is that all there is?} \label{sec:isthatall}

Seen in the context of the other doughnuts, Morley's Theorem looks more
mundane than miraculous.
To bring back the magic,
it would help to fill in the rest of the doughnut holes.
Is there really no way to fill them in?

The cyclotomic identities that underlie Conway's doughnuts
are the basis for the Kubert identities for hyperbolic volume.
(See \cite{milnor:hyperbolic, mohanty:thesis}.)
Dupont and Sah
\cite{dupontsah}
showed that in principle all the Kubert identities
have scissors congruence proofs.
Mohanty
\cite{mohanty:thesis}
worked out explicitly the proofs you get from small $n$.
The results are more complicated than you would hope,
and far removed from anything you would simply guess.

Presumably the algebraic proof of the cyclotomic identities
could be made geometric by turning each step of the algebra into
geometry.
From here we might find that in principle it is always possible to fill
the doughnut hole,
though here again the results might not be as nice as you would hope.

\clearpage
\section{Flip book} \label{sec:flip}
Here are doughnuts.
Download this document and click through the 
following pages to see an animation.
(Or print the pages and flip through them?)
They are arranged in decreasing order of $n$ to suggest that we should
think of Conway's Theorem, Morley's Theorem, and the Bisector Theorem
as special cases of Conway's doughnuts.

\newcommand{\flip}[1]{\newpage\includegraphics[width=370pt]{figures/flip#1.pdf}}
\flip{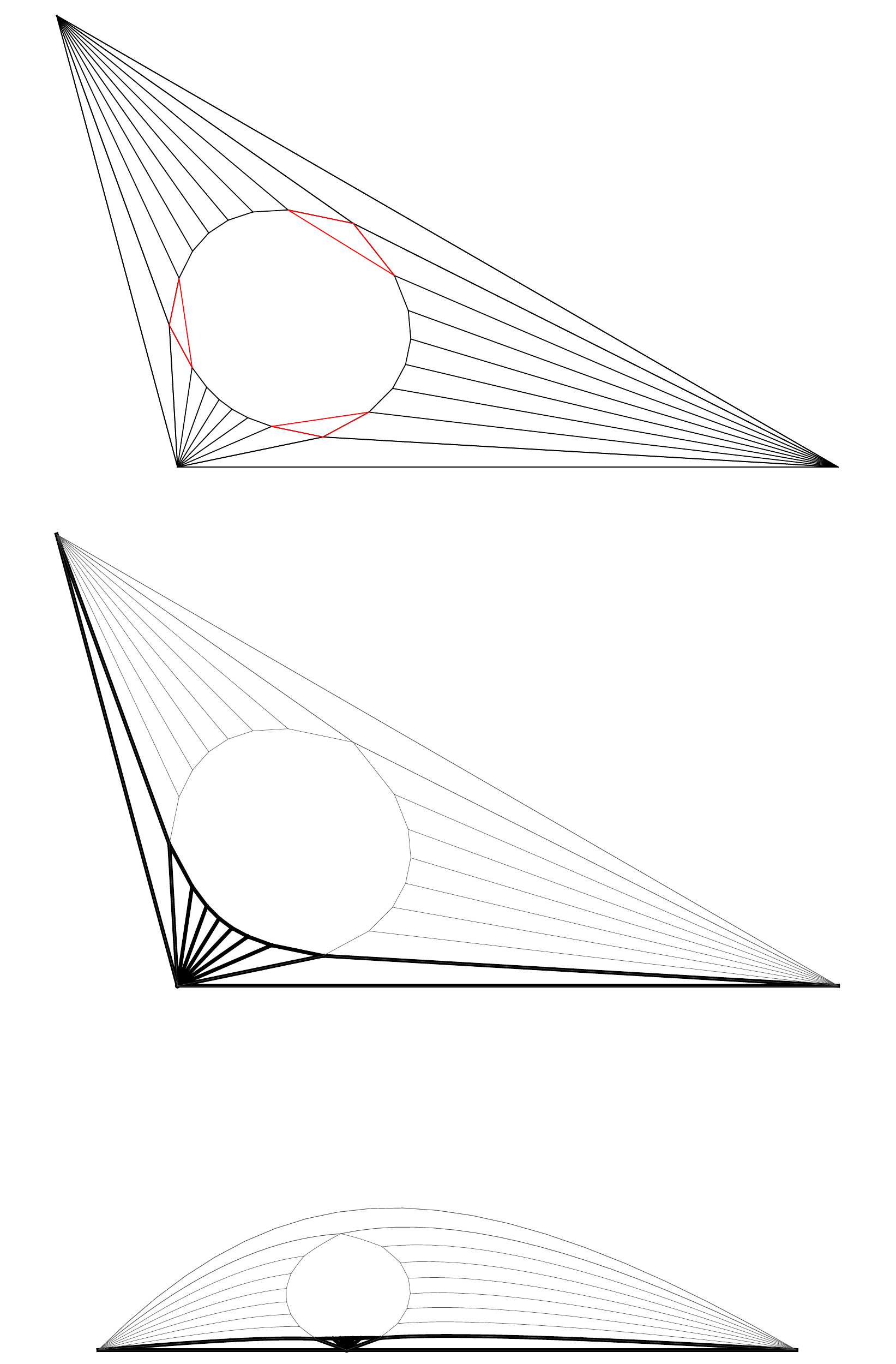}
\flip{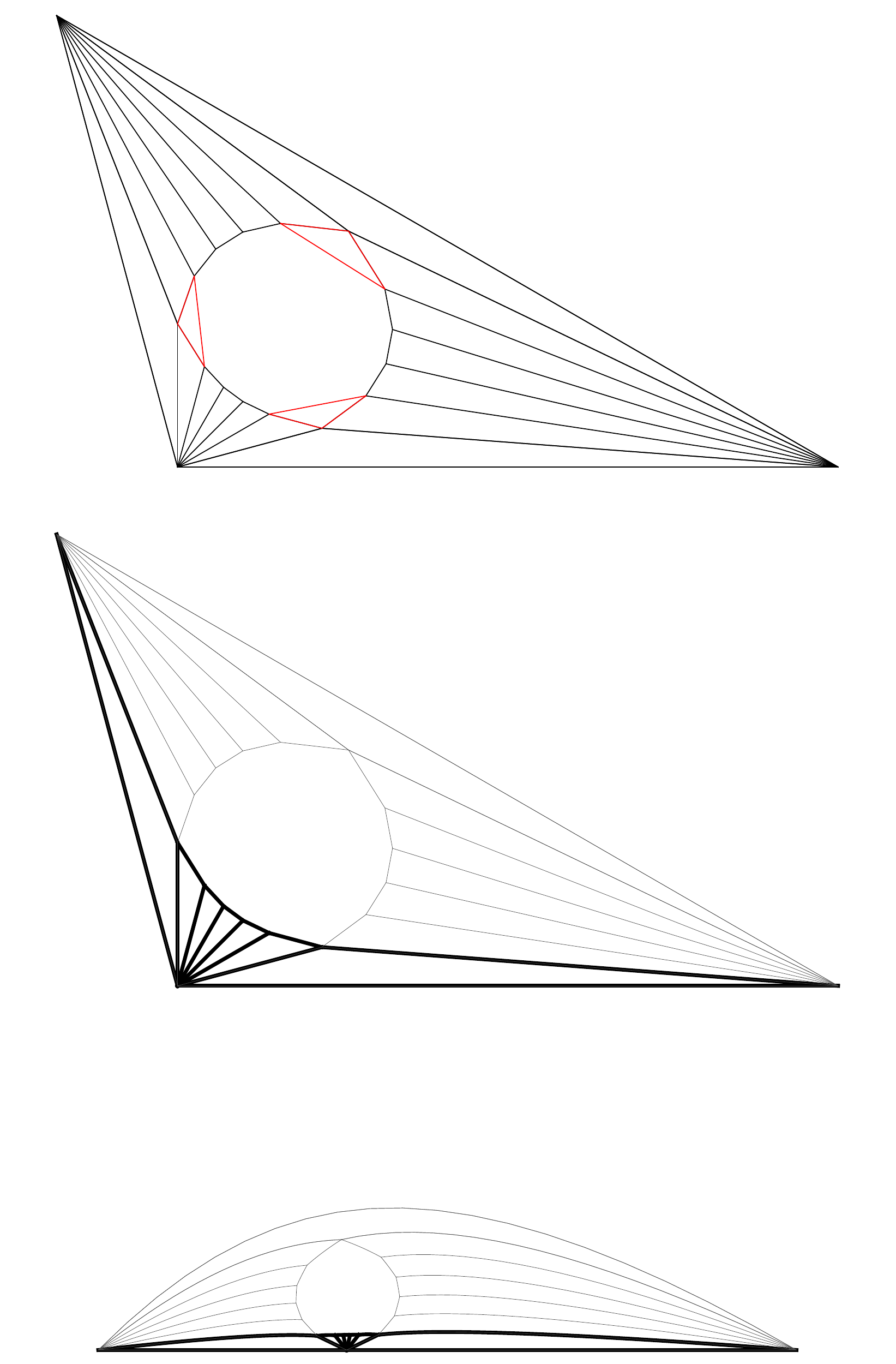}
\flip{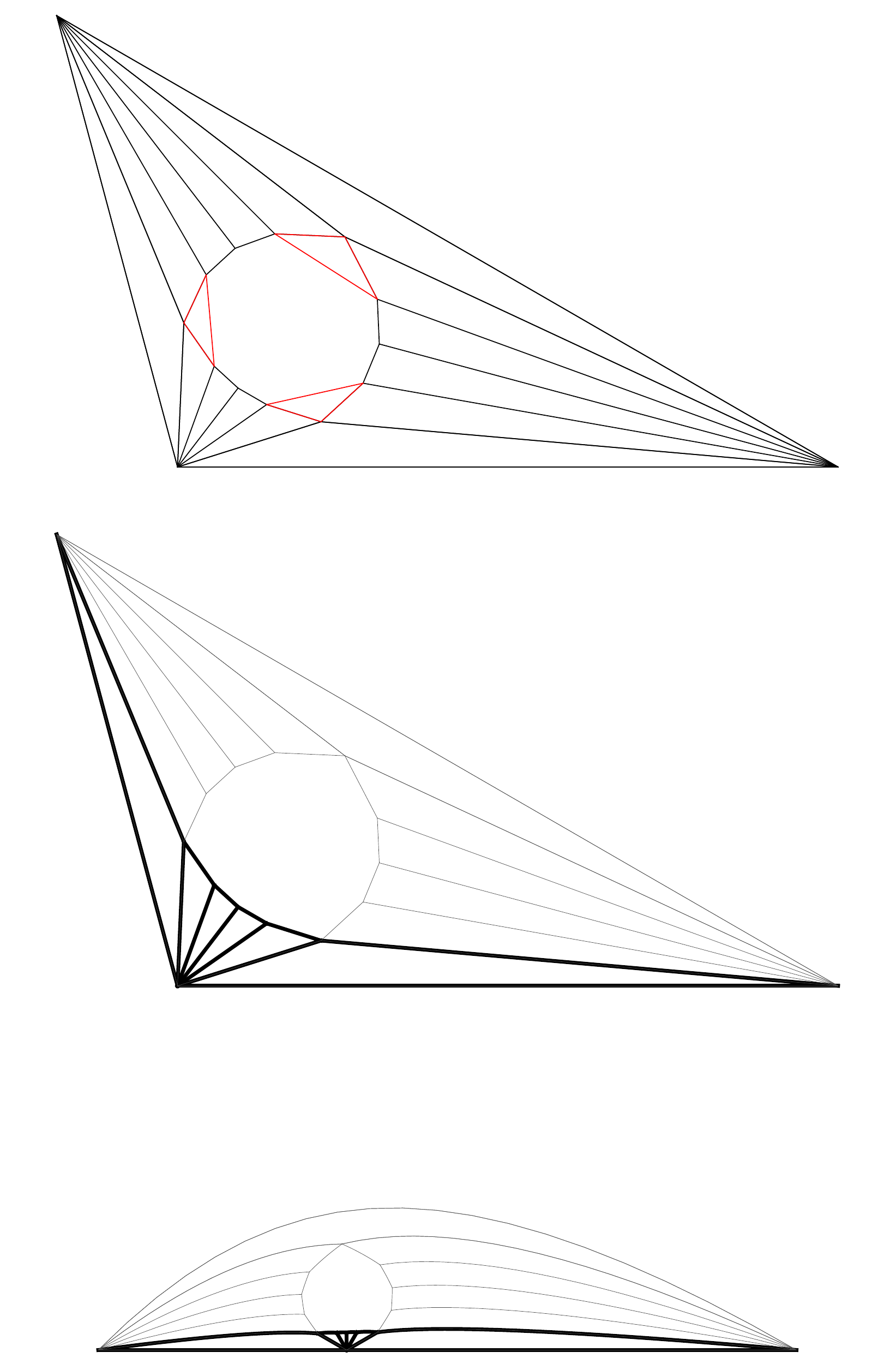}
\flip{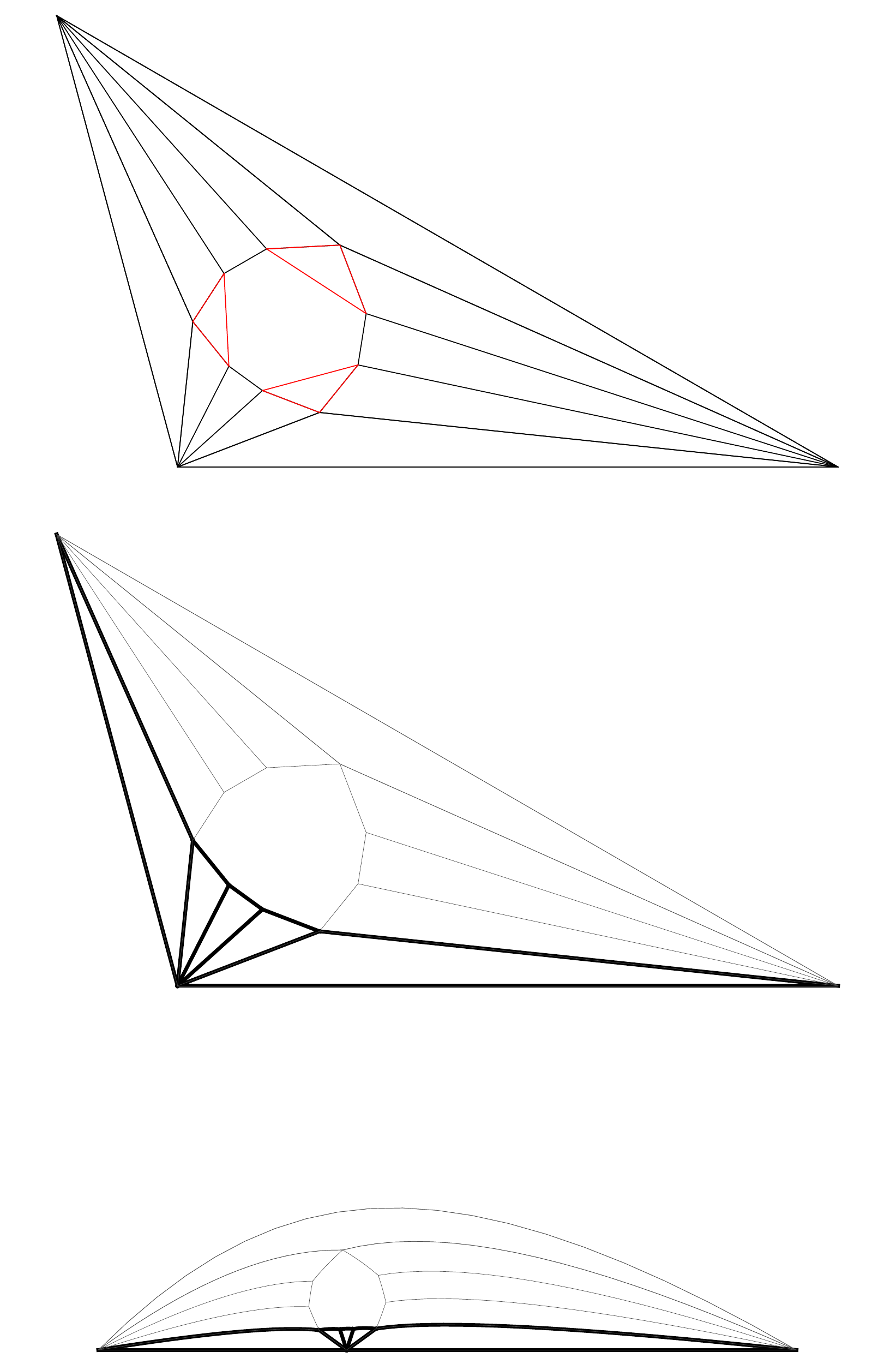}
\flip{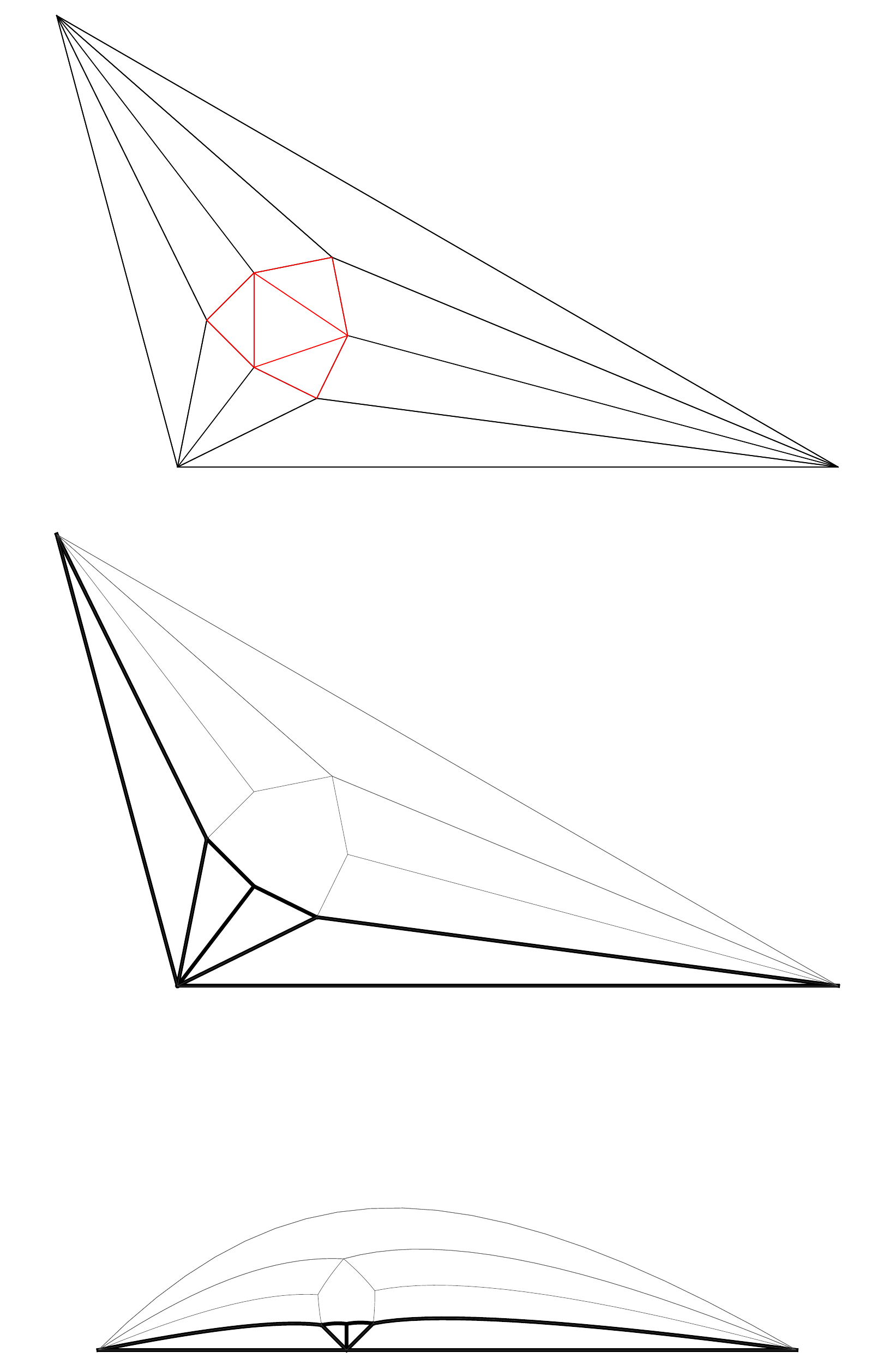}
\flip{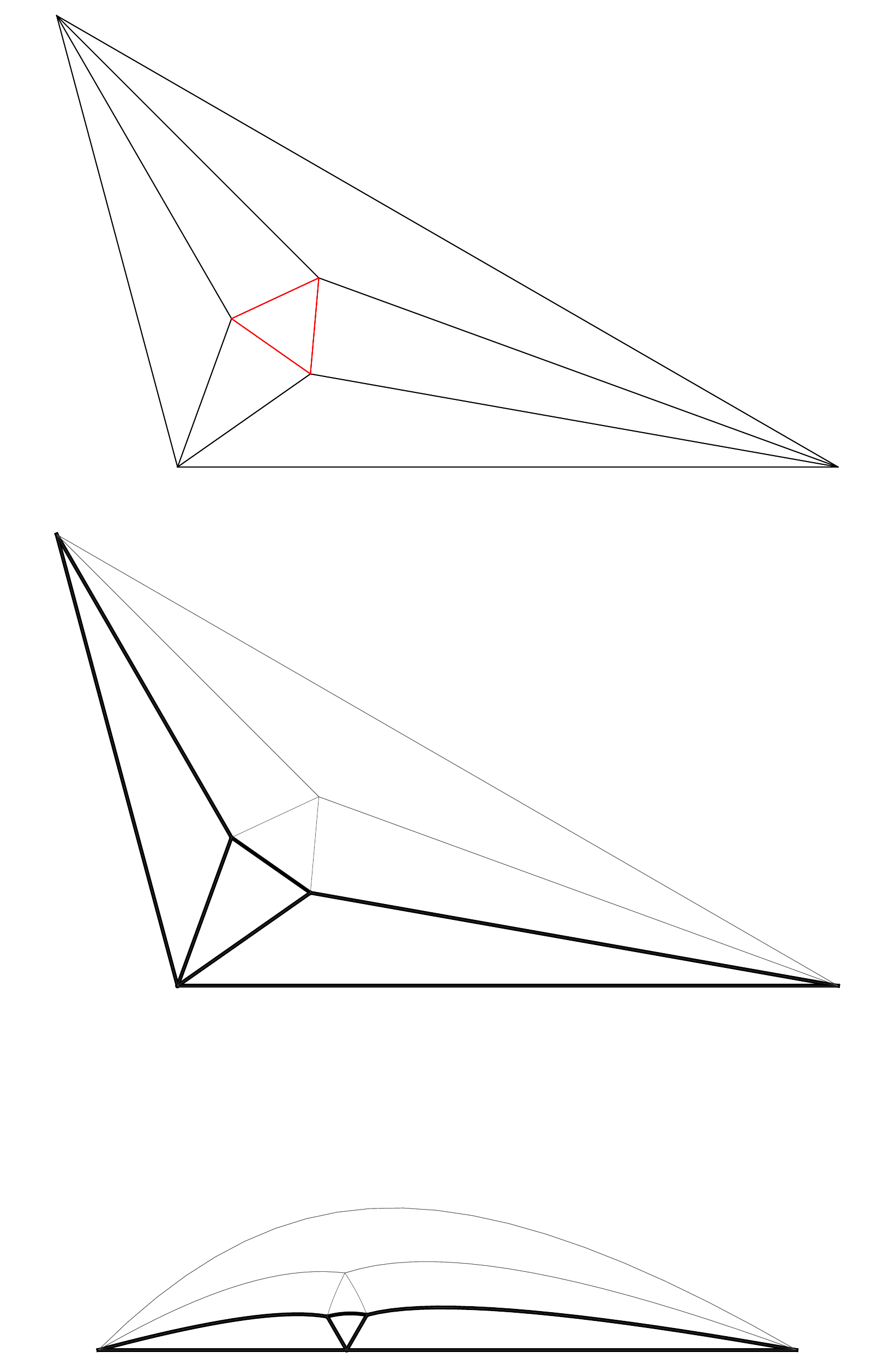}
\flip{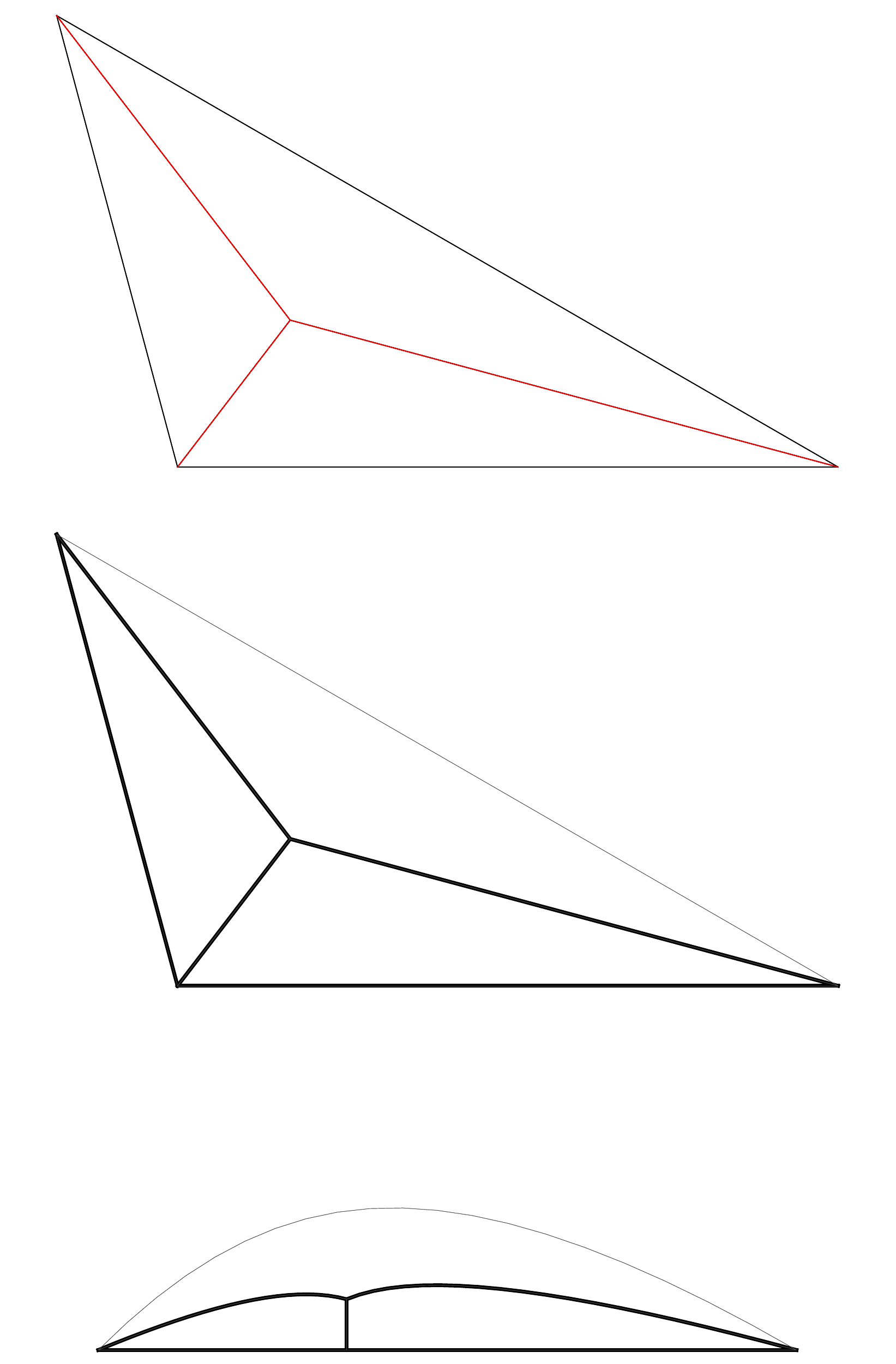}
\flip{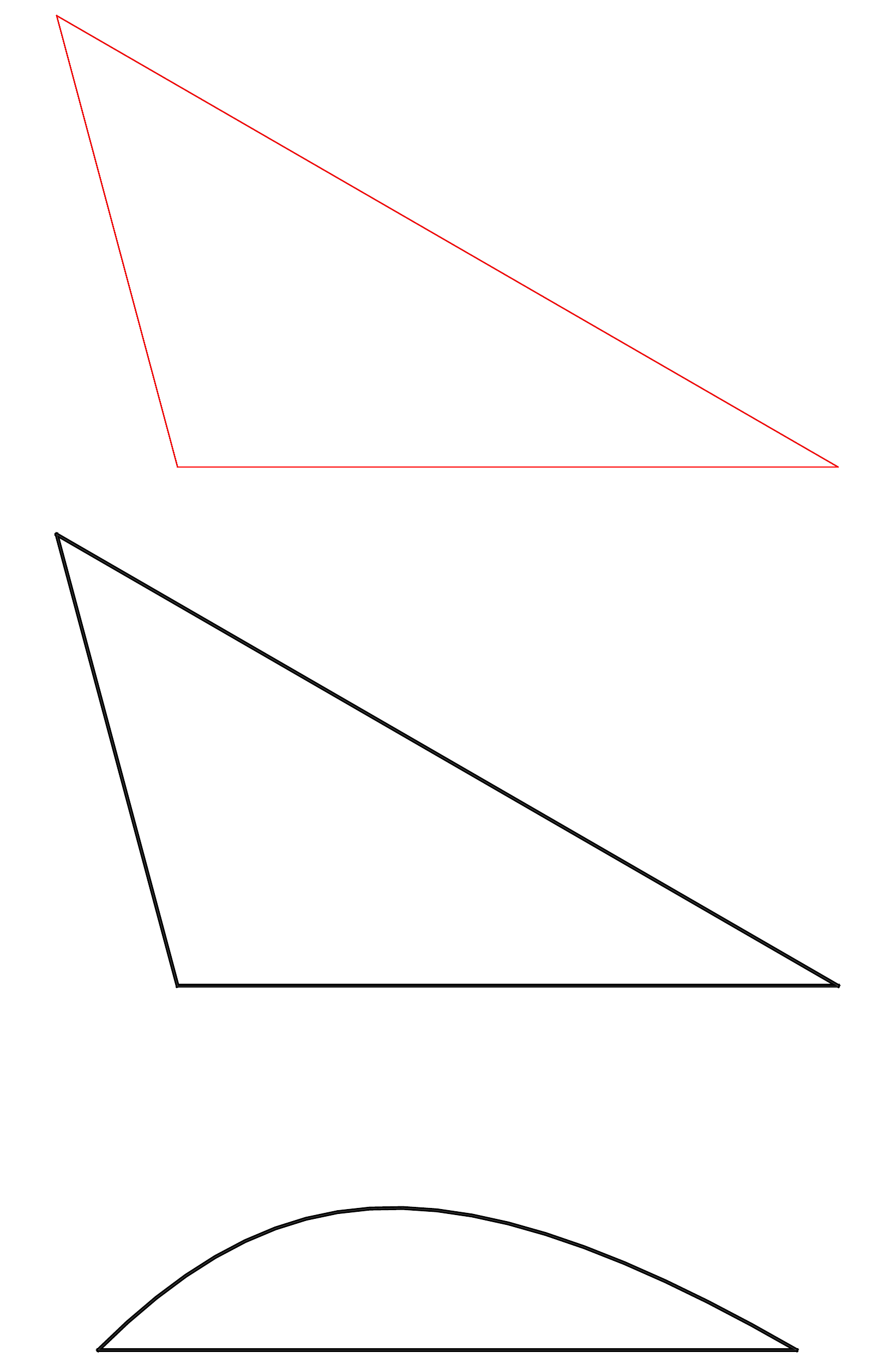}
\clearpage
\bibliographystyle{hplain}
\bibliography{dough}
\end{document}